\documentclass{article}
\usepackage{amsfonts}

\usepackage{amscd}
\usepackage{amsmath}
\usepackage{amssymb}
\usepackage{color}
\usepackage{array}
\usepackage{theorem}
\usepackage{epsfig}
\usepackage{verbatim}

\def\emty{\emptyset}
\def\cbar{\overline{\C}}

\def\stab{\mbox{Stab}\,}

\def\FFF{{\cal F}}

\def\HHH{{\cal H}}

\def\NNN{{\cal N}}

\def\RRR{{\cal R}}

\def\YYY{{\cal Y}}

\def\th{theorem }

\def\Th{Theorem }

\def\polys{polynomials }
\def\qc{quasiconformal }

\def\homeo{homeomorphism }
\def\homeos{homeomorphisms }

\def\cc{connected component }
\def\ccs{connected components }
\def\nbhd{neighbourhood }
\def\nbhds{neighbourhoods }

\def\be{\beta}

\def\g{\gamma}

\def\si{\sigma}

\def\vp{\varphi}

\def\La{\Lambda}
\def\De{\Delta}
\def\de{\delta}

\def\oo{\mbox{$\Omega$} }

\def\R{\mbox{$\mathbb R$}}

\def\C{\mbox{$\mathbb C$}}

\def\Z{\mbox{$\mathbb Z$}}

\def\SS{\mbox{$\mathbb S$}}

\def\DD{\mbox{$\mathbb D$}}

\newtheorem{thm}{Theorem}
\newtheorem{lemma}{Lemma}
\newtheorem{cor}{Corollary}
\newtheorem{prop}{Proposition}
\newcommand{\bx}[1]{\bibitem{#1}}

\title{Pinching Holomorphic Correspondences}
\author{Shaun Bullett\\
School of Mathematical Sciences\\
Queen Mary, University of London\\
Mile End Road, London E1 4NS, UK \\ \\
Peter Ha\"issinsky\\
LATP/CMI\\
Universit\'e de Provence\\
39 rue Fr\'ed\'eric Joliot-Curie\\
13453 Marseille Cedex 13, France}

\begin{document}

\maketitle

\begin{abstract}
For certain classes of holomorphic correspondences which are
matings between Kleinian groups and polynomials, we prove the
existence of pinching deformations, analogous to Maskit's
deformations of Kleinian groups which pinch loxodromic elements to
parabolic elements. We apply our results to establish the
existence of matings between quadratic maps and the modular group,
for a large class of quadratic maps, and of matings between the
quadratic map $z\to z^2$ and circle-packing representations of the
free product $C_2*C_3$ of cyclic groups of order $2$ and $3$.
\end{abstract}



\tableofcontents

\section{Introduction}

It is a well-known consequence of the simultaneous uniformisation
theorem of Bers \cite{Bers} that given two abstractly isomorphic
Fuchsian groups $G_1\subset PSL_2(\R)$ and $G_2\subset PSL_2(\R)$,
acting on the upper and lower complex half-planes respectively,
each having limit set $\hat{\mathbb R}={\mathbb R} \cup \infty$,
and such that the action of $G_1$ on $\hat{\mathbb R}$ is
topologically conjugate to that of $G_2$, the actions of $G_1$ and
$G_2$ can be {\it mated} to obtain a quasifuchsian Kleinian group
$G \subset PSL_2(\C)$. This {\it mating} is a group which is
abstractly isomorphic to both $G_1$ and $G_2$, it has limit set
$\Lambda(G)$ a simple closed (fractal) curve, and the actions of
$G$ on the two components of $\Omega={\hat{\mathbb C}}-\Lambda$
are conformally conjugate to those of $G_1$ on ${\mathcal U}$ and
$G_2$ on ${\mathcal L}$.

\medskip
It is also well-known that given two polynomial maps $P$ and $Q$
of the same degree $n$, in appropriate circumstances one can find
a rational map $R$ which realises a {\it mating} between the
actions of $P$ and $Q$ on their filled Julia sets, in a precise
sense as defined for example in \cite{HT}. A necessary condition
for a mating between two quadratic polynomials $P:z \to z^2+c$ and
$Q:z \to z^2+c'$ to exist is that $c$ and $c'$ should not belong
to conjugate limbs of the connectivity locus (the Mandelbrot Set)
in parameter space: this was first shown also to be a sufficient
condition in the case that $P$ and $Q$ are {\it postcritically
finite} \cite{ree,T}, and subsequently for much more general
classes of $P$ and $Q$ \cite{HT}.

\medskip
In \cite{BP} the first examples of holomorphic correspondences
realising {\it matings} between Fuchsian groups and polynomials
were presented. {\it Holomorphic correspondences} on the Riemann
sphere are multi-valued maps $f:z\to w$ defined by polynomial
equations $p(z,w)=0$. Examples of holomorphic correspondences are
those defined by a union of the graphs of some finite set of
M\"obius transformations, or by the graph of a rational map (or
its inverse). We say that such a correspondence has {\it bidegree}
$(m:n)$ if a generic point $z$ has $n$ images $w$ and a generic
point $w$ has $m$ inverse images $z$.

\medskip
{\bf Definition} {\it Let $q_c:z \to z^2+c$ be a quadratic map
with connected filled Julia set $K(q_c)$. A holomorphic
correspondence $f:z \to w$ of bidegree $(2:2)$ is called a {\it
mating} between $q_c$ and the modular group $PSL_2(\Z)$ if:

\medskip
(a) there exists a completely invariant open simply-connected
region $\Omega \subset \hat{\C}$ and a conformal bijection $h$
from $\Omega$ to the upper half-plane conjugating the two branches
of $f\vert_\Omega$ to the pair of generators $z \to z+1,\ z \to
z/(z+1)$ of $PSL_2(\Z)$;

\medskip
(b) the complement of $\Omega$ is the union of two closed sets
$\Lambda_-$ and $\Lambda_+$, which intersect in a single point and
are equipped with homeomorphisms $h_\pm: \Lambda_\pm \to K(q_c)$,
conformal on interiors, respectively conjugating $f$ restricted to
$\Lambda_-$ as domain and codomain to $q_c$ on $K(q_c)$, and
conjugating $f$ restricted to $\Lambda_+$ as domain and codomain
to $q_c^{-1}$ on $K(q_c)$.}

\medskip
In \cite{BP} the one parameter family of correspondences

$$\left(\frac{az+1}{z+1}\right)^2+\left(\frac{az+1}{z+1}\right)
\left(\frac{aw-1}{w-1}\right)+\left(\frac{aw-1}{w-1}\right)^2=3
\leqno(1) $$

was shown to contain examples of matings between quadratic maps
and the modular group. The following conjecture is implicit in the
discussion in Sections 1 and 6 of that paper.

\medskip
{\bf Conjecture 1} {\it The family $(1)$ of $(2:2)$
correspondences contains matings between $PSL_2(\Z)$ and {\it
every} quadratic polynomial having a connected Julia set, that is
to say every $z \to z^2+c$ with $c \in {\mathcal M}$, the {\it
Mandelbrot set}.}

\medskip
Supporting evidence was provided by proofs for particular examples
and numerical experiments suggesting the resemblance between the
space of matings and the Mandelbrot set. However difficulties in
adapting the theory of {\it polynomial-like maps} \cite{DH} to the
setting of {\it pinched polynomial-like maps} prevented a proof.

\medskip
A different question turned out to be easier to answer. The
modular group may be considered as a representation of the free
product $C_2*C_3$ of cyclic groups, of orders two and three, in
$PSL_2(\C)$. Up to conjugacy there is a one parameter family of
such representations and in the parameter space there is a set
${\mathcal D}$, homeomorphic to a once-punctured closed disc, for
which the representation is discrete and faithful. The modular
group corresponds to a particular {\it boundary} point of
$\mathcal D$. Let $r$ be any representation of $C_2*C_3$
corresponding to a parameter value in the {\it interior}
${\mathcal D}^\circ$ of $\mathcal D$. The ordinary set $\Omega(r)$
of the Kleinian group defined by such a representation $r$ is
connected and the limit set $\Lambda(r)$ is a Cantor set. In
\cite{BH} the notion of a mating between such a representation $r$
of $C_2*C_3$ and a quadratic polynomial $q_c: z \to z^2+c$ was
introduced: $\Lambda_-$ and $\Lambda_+$ are now disjoint, and
their complement $\Omega$ is canonically associated to $\Omega(r)$
(see Section 2.2). By the application of
polynomial-like mapping theory the following analogue of
Conjecture 1 was proved in \cite{BH}.

\medskip
\begin{thm}\label{mating} For every quadratic map $q_c:z \to
z^2 + c$ with $c \in {\mathcal M}$ and every faithful discrete
representation $r$ of $C_2*C_3$ in $PSL_2(\C)$ having connected
ordinary set, there exists a polynomial relation $p(z,w)=0$
defining a $(2:2)$ correspondence which is a mating between $q_c$
and $r$.
\end{thm}

An outline of the proof of \Th \ref{mating} is presented in
Section 2.2, as a prelude to applying pinching techniques to
the matings it shows to exist.

\medskip
We describe an involution $J$ on $\hat{\mathbb C}$ as {\it
compatible} with a mating $f$ if $(J\circ f) \cup I_{\hat{\mathbb
C}}$ is an equivalence relation, where $I_{\hat{\mathbb C}}$
denotes the identity map on $\hat{\mathbb C}$ and $(J\circ f) \cup
I_{\hat{\mathbb C}}$ denotes the $3:3$ correspondence defined by
the algebraic curve
$$p(z,J(w))(z-w)=0$$ (Here $p(z,w)=0$ is the curve defining $f$.)

\begin{prop}\label{compatible} Every mating with a compatible
involution is conjugate to a correspondence in the following two
parameter family (also considered in \cite{BP}):
$$\left(\frac{az+1}{z+1}\right)^2+\left(\frac{az+1}{z+1}\right)
\left(\frac{aw-1}{w-1}\right)+\left(\frac{aw-1}{w-1}\right)^2=3k
\leqno(2) $$
\end{prop}

\medskip
As we shall see, the matings constructed in \cite{BH} have
compatible involutions, so they have representatives in the family
(2), a fact observed in \cite{BH} but for which Proposition
\ref{compatible} (proved in Section 2) provides a more conceptual
setting.

\medskip
The basic idea of pinching can be seen in the process by which the
modular group can be obtained from any chosen standard
representation $r_*$ of $C_2*C_3$ lying in the interior of
$\mathcal D$, that is to say a faithful discrete representation
with connected ordinary set $\Omega(r_*)$ (and therefore limit set
a Cantor set). We first recall that each Kleinian representation
of $C_2*C_3$ comes equipped with a canonical involution $\chi$
which conjugates the generators $\sigma\in C_2$ and $\rho\in C_3$
to their inverses (see Section 2.1); we let $G$ denote the group
$<\chi,\sigma,\rho>$. For each rational number $p/q$ there is an
arc $\delta_{p/q}$ on the orbit space $\Sigma=\Omega(r_*)/G$ which
lifts to simple closed geodesic $\tilde{\delta}_{p/q}$ of winding
number $p/q$ on a certain torus $\tilde{\Sigma}$ double-covering
$\Sigma$ (see Lemma \ref{arcs-exist} in Section 3.1 for details).
The arc $\delta_{p/q}$ lifts to an arc $\alpha_{p/q}$ on
$\Omega(r_*)$ together with its translates under $G$. This arc
$\alpha_{p/q}$ is {\it precisely $<<g>>$-invariant} for any
loxodromic element $g \in G$ which stabilises it. (Here $<<g>>$
denotes the maximal elementary subgroup of $G$ containing $g$, and
saying that an arc $\alpha$ is {\it precisely $<<g>>$-invariant}
means that $<<g>>\alpha=\alpha$ and $h(\alpha)\cap \alpha =
\emptyset$ for all $h\in G$ not in $<<g>>$). In this situation
Maskit's Theorem \cite{M} states that the representation of $G$ in
$PSL_2(\C)$ can be deformed to one in which $\alpha_{p/q}$ and its
translates under $G$ are pinched to points and $g$ becomes
parabolic. We deduce that we may pinch $\delta_0$, and hence its
lift $\alpha_0$, to a point, thereby deforming the representation
$r_*$ of $C_2*C_3$ to the representation $PSL_2(\Z)$, which lies
on the boundary of the deformation space $\mathcal D$. Similarly
for $p/q \ne 0$ we may pinch $\delta_{p/q}$ to a point and so
deform the representation $r_*$ to a faithful discrete
representation which we denote $r_{p/2q}$. This has ordinary set a
disjoint union of a countable infinity of open round discs, and
limit set a circle-packing. The representation $r_{p/2q}$ depends
only on the value of $p/2q$ mod $2$: pinching $\delta_{(2nq+p)/q}$
in place of $\delta_{p/q}$ amounts to approaching the same limit
representation $r_{p/2q}$ but by a non-isotopic path in ${\mathcal
D}$. We remark that by a deep result of McMullen \cite{mc2} the
representations $r_{p/2q}$ are dense in the boundary of ${\mathcal
D}$.

\medskip
Recently, Ha\"{\i}ssinsky \cite{H2}, Cui \cite{cui} and
Ha\"{\i}ssinsky and Tan \cite{HT} proved  analogous results to
Maskit's in the context of rational maps, showing that, under
appropriate hypotheses, given a rational map $R$ and an
$R$-invariant union of arcs joining attracting to repelling
cycles, one can continuously deform the map in such a way that the
arcs, and their pre-images, are pinched to points and the cycles
become parabolic.

\medskip
In this paper we adapt the techniques of \cite{H2} and \cite{HT}
to apply them to the holomorphic correspondences constructed in
\cite{BH}. In Section 3, for any correspondence $p_0(z,w)=0$ which
is a mating between $r_*$ and $q_c$, and for any rational number
$p/q$, we identify an arc $\gamma_{p/q}$ such that the grand orbit
of $\gamma_{p/q}$ under the correspondence is a union of
infinitely many disjoint copies of $\gamma_{p/q}$ (or of copies of
a quotient of $\gamma_{p/q}$ by an involution), and such that
pinching each connected component of this union to a point
corresponds to deforming the representation $r_*$ to $r_{p/2q}$.
We describe the pinching process formally as follows.

\medskip
{\bf Definition} {\it A convergent pinching deformation for
$\gamma_{p/q}$ is a family of \qc maps $(\vp_t)_{0\le t <1}$ of
the Riemann sphere such that the conjugate correspondences $p_t$
defined by
$$p_t(z,w)= p_0(\vp_t^{-1}(z),\vp_t^{-1}(w))$$
are holomorphic and satisfy the following\,:
\begin{itemize}
\item $(p_t,\vp_t)$ are uniformly convergent to a pair
$(p_1,\vp_1)$ as $t$ tends to $1$ , \item the non-trivial fibres
of $\vp_1$ are exactly the closure of the \ccs of the orbit of
$\g_{p/q}$. \end{itemize}}

\medskip
There are two technical conditions that we require the quadratic
map $q_c$ to satisfy in order to apply the techniques of \cite{HT}
to $\gamma_0$:

\medskip
(i) if the critical point $0$ of $q_c$ is recurrent, the
$\beta$-fixed point of $q_c$ is not in the $\omega$-limit set of
$0$;

\medskip
(ii) $q_c$ is weakly hyperbolic, that is, there are constants
$r>0$ and $\de<\infty$ such that, for all $z\in J_q\smallsetminus
\{\text{preparabolic points}\}$, there is a subsequence of
iterates $(q^{n_k})_k$ such that
$$\mbox{deg}(W_k(z)\stackrel{q^{n_k}}{\longrightarrow} D(q^{n_k}(z),r) )\le
 \de\,$$
where $W_k(z)$ is the connected component of
$q^{-n_k}(D(q^{n_k}(z),r) )$ containing $z$.

\medskip
In Section 4 we prove:

\medskip
\begin{thm}\label{simple}
Let $p_0(z,w)$ be a mating between the representation $r_*$ and
$q_c$, where $q_c$ satisfies conditions (i) and (ii) above. Then
there exists a pinching deformation of $p_0$ such that $(p_t)_{0
\le t <1}$ converges uniformly to a mating $p_1$ between
$PSL_2(\Z)$ and $q_c$.\end{thm}

\begin{figure}
\begin{center}
\includegraphics{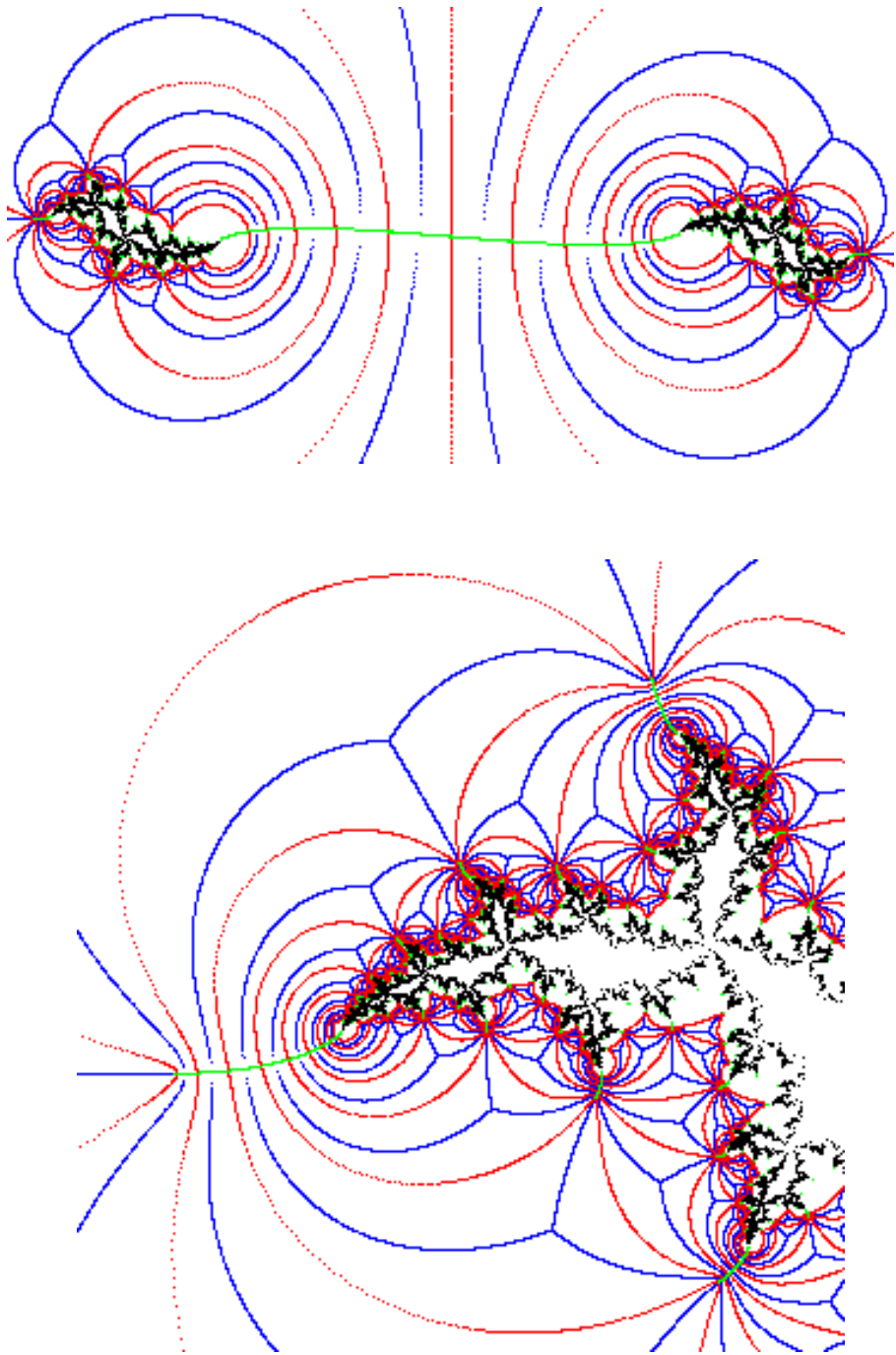}
\caption{A mating of a representation of $C_2*C_3$ with a Douady
rabbit (and zoom). The arc $\gamma_0$ and its images are shown.
Pinching these gives a mating of $PSL_2(\Z)$ with the rabbit, by
Theorem 2.}
\end{center}
\end{figure}

\medskip
{\bf Corollary} {\it Conjecture 1 is true for all quadratic maps
$q_c$ which satisfy conditions (i) and (ii).}

\medskip
The class of {\it weakly hyperbolic} quadratic maps is quite
large: for example it contains all quadratic maps which satisfy
the {\it Collet-Eckmann} condition \cite{pr}, and those which
contain parabolic points.

\medskip
We next investigate pinching $\gamma_{p/q}$, for $p/q \ne 0$. In
Section 3.2, we define the notion of a  mating between the
circle-packing representation  $r_{p/2q}$ of $C_2*C_3$ and $q_c$.
This generalises our earlier definition of a {\it mating between
$PSL_2(\Z)$ and $q_c$}, replacing $K(q_c)$ by a certain
identification space $K(q_c)/\sim_{p/q}$ and replacing the
condition that $\Lambda_+\cap\Lambda_-$ be a point by the
condition that it consist of $q$ points (the $p/q$ {\it Sturmian
orbit} on the boundary of $K(q_c)$). We show that a mating between
$r_{p/2q}$ and $q_c$ depends only on $p/q$ mod $1$. To avoid
technical difficulties we restrict attention to the special case
that the quadratic map is $q_0:z \to z^2$. Using the techniques of
\cite{H2}, we prove the following:

\begin{thm}\label{rat} Let $p_0(z,w)$ be a mating between the
representation $r_*$ and $q_0$, and let $p/q$ be any rational.
Then there exists a pinching deformation of $p_0$ such that
$(p_t)_{0 \le t <1}$ converges uniformly to a mating $p_1$ between
the circle-packing representation $r_{p/2q}$ of $C_2*C_3$ in
$PSL_2(\C)$ and $q_0$.
\end{thm}

The following is the natural generalisation of Conjecture 1.

\medskip
{\bf Conjecture 2} {\it For every $0\le p/q<1$, the family $(2)$
of $(2:2)$ correspondences contains matings between the
circle-packing representation $r_{p/2q}$ and every quadratic
polynomial $q_c$ which has $c \in {\mathcal M}\setminus{\mathcal
M}_{1-p/q}$, where ${\mathcal M}_{1-p/q}$ denotes the
$(1-p/q)$-limb of the Mandelbrot set ${\mathcal M}$.}

\medskip
The condition that $c$ does not lie in ${\mathcal M}_{1-p/q}$ is
necessary for elementary topological reasons. One might hope to
generalise the techniques of the present paper to prove Conjecture
2 in the case that $q_c$ satisfies conditions (i) and (ii) of the
hypotheses of \Th \ref{simple}, but the technical details could be
formidable.

\medskip
{\bf Warning} As will already be apparent, certain of the
constructions and results in this article depend on $p/q \in
{\mathbb Q}$, certain depend only on $p/q$ mod $1$ (the class of
$p/q$ in ${\mathbb Q}/{\mathbb Z}$), and certain on $p/q$ mod $2$.
We shall try to make the dependence clear in each case, but
briefly the situation may be summed up as follows. A
circle-packing representation $r_{p/2q}$ of $C_2*C_3$ depends on
$p/q$ mod $2$ but the route to it (in the moduli space ${\mathcal
D}$) given by pinching $\delta_{p/q}$ depends on $p/q \in {\mathbb
Q}$. A mating between $r_{p/2q}$ of $C_2*C_3$ and $q_c$ depends
only on $p/q$ mod $1$, but again the route to it (in mating space)
given by pinching $\gamma_{p/q}$ depends on $p/q \in {\mathbb Q}$.

\section{Matings between quadratic maps and representations of $C_2*C_3$}

We define what we mean by {\it matings} between quadratic maps and
representations of $C_2*C_3$ in $PSL_2(\C)$ which lie in
${\mathcal D}^o$, we recall the main ideas of the proof \cite{BH}
of \Th \ref{mating}, we prove Proposition \ref{compatible}, and we
present a group-theoretic description of the `ordinary set'
$\Omega(f)$ of a mating.

\subsection{Faithful discrete
representations with connected ordinary sets}

Up to conjugacy each representation $r$ of $C_2*C_3$ in
$PSL_2(\C)$ is determined by a single complex parameter, the
cross-ratio between the fixed points on $\hat{\mathbb C}$ of the
action of the generator $\sigma$ of $C_2$ and those of the
generator $\rho$ of $C_3$. Such a representation comes equipped
with a (unique) involution $\chi$ which exchanges the two fixed
points of $\sigma$ and also those of $\rho$, so that
$\chi\sigma=\sigma\chi$ and $\chi\rho=\rho^{-1}\chi$. On the
Poincar\'e $3$-ball $\chi$ is simply rotation through $\pi$ around
the common perpendicular to the axes of $\sigma$ and $\rho$. Write
$G$ for the group $<\sigma,\rho,\chi>$, and note that it has
ordinary set $\Omega(G)$ the same as that of $<\sigma,\rho>$.

\medskip The faithful discrete actions $r:C_2*C_3 \subset
PSL_2(\C)$ with connected ordinary set $\Omega(G)$ form a single
quasiconformal conjugacy class, the class of representations for
which one can find simply-connected fundamental domains for
$\sigma$ and $\rho$ with interiors together covering the whole
Riemann sphere (the conditions of the simplest form of the Klein
Combination Theorem are satisfied) \cite{mas}. Such fundamental
domains may be constructed as illustrated in figure 2.

\begin{figure}
\begin{center}
\input{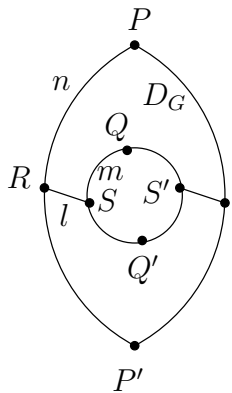}
\caption{A fundamental domain $D_G$ for the group
$G=<\sigma,\rho,\chi>$.} \label{fund}
\end{center}
\end{figure}

\medskip
Here $P$ and $P'$ are the fixed points of $\rho$, $Q$ and $Q'$ are
the fixed points of $\sigma$, $R$ is a fixed point of (the
involution) $\chi\rho$ and $S$ and $S'$ are the fixed points of
$\chi\sigma$. The lines $l,m$ and $n$, joining $R$ to $S$, $Q$ to
$S$ and $R$ to $P$, are chosen such that they are smooth and
remain non-intersecting in the quotient orbifold $\Omega(G)/G$.
The region bounded by $n,\rho n, \chi n$ and $\chi\rho n$ is a
fundamental domain for $\rho$, and the region exterior to the loop
made up of $m,\sigma m,\chi m$ and $\chi\sigma m$ is a fundamental
domain for $\sigma$. The intersection of these two regions is a
fundamental domain for the (faithful) action of $C_2*C_3$ on
$\Omega(G)$, and the half $D_G$ of this intersection bounded by
$n,l,m,\sigma m, \chi l$ and $\rho n$ is a fundamental domain for
the action of $G$. The union of all translates of $D_G$ under
elements of $C_2*C_3$ is a topological disc $D$ which is a
fundamental domain for the action of $\chi$ on $\Omega(G)$. The
complement $\Lambda(G)$ of $\Omega(G)= D \cup \chi(D)$ in
$\hat{\mathbb C}$ is a Cantor set.

\medskip
The orbifold $\Omega(G)/G$ is a sphere $\Sigma$, which has a
complex structure with four cone points, which we may also label
$P,Q,R,S$, where $P$ has angle $2\pi/3$ and $Q,R,S$ each have
angle $\pi$. For a given representation of $C_2*C_3$, a set of
lines $l,m,n$ as in figure 2 descend to an isotopy class of
non-intersecting paths joining the corresponding cone points in
$\Sigma$. By considering the choices we may make of the various
labels and lines in figure 2 we can obtain a description of
$\widetilde{\mathcal D}^0$, the universal cover of the moduli
space ${\mathcal D}^0$.

\begin{lemma}\label{Kleinian-markings}
There is a homeomorphism $\Phi$ between ${\mathcal D}^o$ and the
space $\mathcal S$ of spheres $\Sigma$ having a complex structure
with four marked cone points $P,Q,R,S$ where $P$ has angle
$2\pi/3$ and $Q,R,S$ each have angle $\pi$. This homeomorphism
$\Phi$ lifts to a homeomorphism $\tilde{\Phi}$ between
$\widetilde{\mathcal D}^o$ and the space $\tilde{\mathcal S}$ of
spheres $\Sigma \in {\mathcal S}$ marked with an isotopy class of
non-intersecting paths $PR$, $RS$ and $SQ$.
\end{lemma}

{\bf Proof.} For $r\in {\mathcal D}^o$ define $\Phi(r)$ to be the
orbifold $\Omega(G)/G$, where $G=<\sigma,\rho,\chi>$ is the
subgroup of $PSL_2(\C)$ corresponding to the representation $r$.
Clearly $\Phi$ is continuous as ${\mathcal D}^o$ is endowed with
the topology induced by its parametrisation by the cross-ratio
$(Q,Q';P,P')$. To define an inverse to $\Phi$, observe that given
any $\Sigma \in {\mathcal S}$, we may obtain a representation $r$
by regarding $\Sigma$ as a quasiconformal deformation of the
orbifold corresponding to $r_*$, lifting the corresponding ellipse
field to $\hat{\C}$, and applying the Measurable Riemann Mapping
Theorem.

\medskip
To lift $\Phi$ to a homeomorphism $\tilde{\Phi}$ we have to
consider markings. Note that given a representation of $C_2*C_3$
which lies in ${\mathcal D}^o$, there is only one choice for which
of the pair $P,P'$ (in figure 2) to label $P$, namely the fixed
point of $\rho$ around which the rotation is {\it anticlockwise}.
There is also just one choice (up to isotopy) for the arc $n$. The
labels $Q$ and $Q'$ are interchangeable (provided that we also
interchange the labels $S$ and $S'$), but once a choice has been
made for $Q$ the arc $m$ is determined, and even if the labels $Q$
and $Q'$ are exchanged the arc $QS$ in the orbifold $\Sigma$ is
unchanged up to isotopy. This just leaves us a choice of the arc
$l$ in figure 2. We can alter $l$ to wind an extra $n$ times
around the central `hole' for any integer $n$, or $n+1/2$ times if
we switch the labels $Q$ and $Q'$. Changing the winding number of
$l$ corresponds to choosing a different isotopy class of paths
between the points labelled $R$ and $S$ in the orbifold $\Sigma$.
$\square$

\medskip
Let $t_\alpha$ denote the automorphism of $\widetilde{\mathcal
D}^o$ corresponding to turning the internal boundary of figure 2
through an angle $2\pi\alpha$. Note that $t_{1/4}$ moves the pair
of points labelled $Q,Q'$ to the pair labelled $S,S'$ and vice
versa. Let $\iota:{\mathcal D}^o \to {\mathcal D}^o$ denote the
involution obtained by replacing the generating pair
$\{\sigma,\rho\}$ of $C_2*C_3$ by $\{\sigma',\rho \}$, where
$\sigma'=\chi\sigma$. This corresponds to composing the
representation with an outer automorphism of $C_2*C_3$. The
following result is now self-evident.

\begin{lemma}\label{quarter-twist} The automorphism
$t_{1/4}:\widetilde{\mathcal D}^o \to \widetilde{\mathcal D}^o$
covers $\iota:{\mathcal D}^o \to {\mathcal D}^o$, and $t_{1/2}$
generates the group of covering transformations of
$\widetilde{\mathcal D}^o \to {\mathcal D}^o$. $\square$
\end{lemma}

\subsection{Matings between $q_c$ and $r \in {\mathcal D}^o$}

As in the previous subsection, $G$ denotes the group
$<\sigma,\rho,\chi>$.

\medskip
{\bf Definition} {\it A $(2:2)$ holomorphic correspondence $f:z
\to w$ is called a {\it mating} between a faithful discrete
representation $r$ of $C_2*C_3$ in $PSL_2(\C)$ having connected
ordinary set $\Omega(G)$ and a polynomial $q_c:z \to z^2+c$ having
connected filled Julia set $K(q_c)$, if the Riemann sphere
$\hat{\mathbb C}$ is the disjoint union of a connected open set
$\Omega(f)$ and a closed set $\Lambda(f)$ made up of two
components, $\Lambda_+(f)$ and $\Lambda_-(f)$ such that each of
$\Omega(f)$ and $\Lambda(f)$ is completely invariant under $f$
and:

\medskip
(a) the action of $f$ on $\Omega(f)$ is discontinuous and there is
a conformal bijection between the grand orbit space $\Omega(f)/f$
and $\Omega(G)/G$;

\medskip
(b) there is a \qc homeomorphism defined from a \nbhd of
$\Lambda_-(f)$ onto a \nbhd of  $K(q_c)$ in $\C$, which realises a
hybrid equivalence, conjugating $f$ to $q_c$. Similarly there is a
hybrid equivalence between $(f^{-1},\Lambda_+(f))$ and
$(q_c,K(q_c))$, this time conjugating $f^{-1}$ to $q_c$.}

\medskip
(See \cite{DH} for the definition of the term `hybrid
equivalence'.)

\bigskip
The construction of a holomorphic correspondence which realises a
mating between given $q_c$ and $r$ proceeds as follows (see
\cite{BH} for more details).

\medskip
We first associate an annulus $A$ to $q_c:z \to z^2+c$. There is a
holomorphic conjugacy (the B\"ottcher coordinate) from $z \to z^2$
to $q_c$ on a neighbourhood of $\infty$, fixing the point $\infty$
and tangent to the identity map there \cite{DH1}. An {\it
equipotential} for $q_c$ is the image of a circle $\{Re^{2\pi it}:
0 \le t <1\}$ under this conjugacy. It is a smooth Jordan curve
parameterized by {\it external angle} $t$. The region bounded by
such an equipotential is a simply-connected domain $V$, mapped
$2:1$ by $q_c$ onto a larger domain $U \supset \overline{V}$ which
also has boundary an equipotential parametrised by external angle.
Let $A$ denote the annulus $U-V$, and denote its inner and outer
boundaries by $\partial_1 A$ and $\partial_2 A$ respectively. The
map $q_c$ sends $\partial_1 A$ two-to-one onto $\partial_2 A$.
There is an involution $i:z \to -z$ on $V$ sending each $z \in V$
to the other point which has the same image in $U$ under $q_c$,
and there are many choices possible of an orientation-reversing
smooth involution $j$ on $\partial_2 A$, a canonical choice being
given by $t \to 1-t$ on external angles.

\begin{figure}
\begin{center}
\input{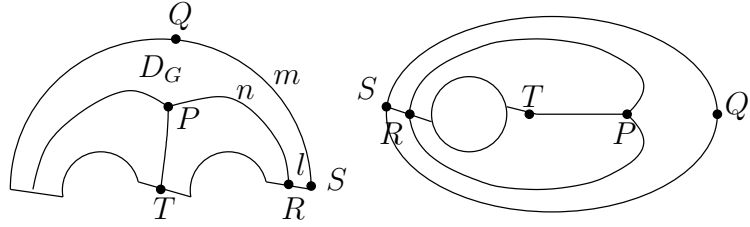}
\caption{The set $D_G \cup \rho(D_G) \cup
\rho^{-1}(D_G)$ and its quotient the annulus $B$.}
\label{ann}
\end{center}
\end{figure}

\medskip
The next ingredient is an annulus $B$ associated to $r$. Recall
the fundamental domain $D_G$ constructed above for the group
$G=<\sigma,\rho,\chi>$ and illustrated in figure 2. Let $B$ denote
the annulus consisting of the three copies $D_G \cup \rho D_G \cup
\rho^{-1}D_G$ of $D_G$, with the boundary identifications (induced
by $\chi$) indicated in figure 3. The rotations $\rho$ and
$\rho^{-1}$ mapping $D_G \cup \rho D_G \cup \rho^{-1} D_G$ to
itself descend to a $2:2$ correspondence $g$ on $B$, mapping each
$z \in B$ to the pair $\{\rho z, \rho^{-1}z\}$ (or rather to their
equivalence classes under the action of $\chi$). The set $D_G$
descends to a `fundamental domain' for the action of $g$ on $B$.
The boundary of $B$ is divided into three segments (two inner and
one outer, figure 3), each of which is mapped to the other two by
$g$. Thus when its domain is restricted to the inner boundary
$\partial_1 B$, and its range is restricted to the outer boundary
$\partial_2 B$, the correspondence $g$ defines a two-to-one map.
When restricted to a correspondence from the inner boundary to
itself, $g$ defines a (fixed point free) bijection. Moreover the
involution $\sigma$ descends to an involution (which we also
denote $\sigma$) on the outer boundary $\partial_2 B$ of $B$,
having fixed points $Q$ and $S$.

\medskip
\begin{lemma}\label{AB} There exists a quasiconformal homeomorphism
$h$ from $A$ to $B$ which restricts to a smooth homeomorphism from
$\partial A$ to $\partial B$ conjugating the boundary maps
$(q_c:\partial_1 A \to
\partial_2 A,\ j:\partial_2 A \to \partial_2 A)$ to the
boundary maps $(\sigma \circ g:\partial_1 B \to
\partial_2 B,\ \sigma: \partial_2 B \to \partial_2 B)$.
\end{lemma}

\medskip This lemma is proved \cite{BH} by applying standard
techniques of Ahlfors and Bers. Now to construct a mating between
$q_c$ and $r$ first glue together $U$ and a second copy $U'$ of
$U$, via the boundary involution $j$, to obtain a sphere $U \cup
U'$, equipped with an involution (also denoted $j$) exchanging $U$
with $U'$ and restricting to the original $j$ on the common
boundary. Inside $U'$ is a simply-connected subdomain $V'$
corresponding to $V \subset U$. Let $q_c'=j\circ q_c \circ j:V'
\to U'$ denote the quadratic map corresponding to $q_c:V \to U$
and $A'$ denote the annulus $U'-V'$. To define a $2:2$ topological
correspondence $f$ on $U \cup U'$ we fit together:

$\bullet$ $q_c:V \to U$ (a $2:1$ correspondence);

$\bullet$ $(q_c')^{-1}=j\circ q_c^{-1} \circ j:U' \to V'$ (a $1:2$
correspondence);

$\bullet$ $j\circ i:V \to V'$ (a $1:1$ correspondence), and

$\bullet$ $j\circ g:A \to A'$ (a $2:2$ correspondence),

where $g:A \to A$ is the $2:2$ correspondence constructed earlier.
Now define an ellipse field on $A$ by using Lemma \ref{AB} to
transport the standard complex structure from the annulus $B$.
Using $j$ extend this ellipse field to $A'$ and pulling back via
$q_c^{-1}$ and $q_c'^{-1}$ extend it to an ellipse field on the
whole of $\hat{\mathbb C}-(K(q_c) \cup K(q_c'))$, which transforms
equivariantly under the action of the $2:2$ correspondence $f$.
Extend this ellipse field to the whole of $\hat{\mathbb C}$ by
using the standard complex structure on $K(q_c) \cup K(q_c')$. By
applying the Measurable Riemann Mapping Theorem we obtain a
complex structure respected by $f$, completing our outline proof
of Theorem \ref{mating}.

\medskip
For any mating $f$ constructed by the method of the proof above,
the $3:3$ correspondence $(j\circ f)\cup I_{\hat{\mathbb C}}$
sends each $z\in V$ to the triple of points $\{z,i(z),jq_c(z)\}$,
each $z\in A$ to the triple $\{z,g(z)\}$ (recall that $g$ is $2:2$
so $g(z)$ contains two points), and each $z \in U'$ to the triple
$\{z,q_c^{-1}j(z)\}$. It is easily checked that each of these
triples is the grand orbit under $(j\circ f)\cup I_{\hat{\mathbb
C}}$ of any one of its elements, in other words the $3:3$
correspondence is an equivalence relation. The involution $j$ is
therefore {\it compatible} with the mating $f$ in the sense
defined in Section 1. To show that $f$ is conjugate to a
correspondence in the family $(2)$ it now only remains to prove
Proposition \ref{compatible}. But a holomorphic correspondence
which is an equivalence relation is necessarily the covering
correspondence of a rational map, and so there is a rational map
$Q$ of degree three such that $(J\circ f)\cup I_{\hat{\mathbb
C}}=Cov^Q$ where
$$Cov^Q: z \to w \quad \Leftrightarrow \quad Q(w)-Q(z)=0.$$
We deduce that $f=J\circ Cov^Q_0$, where
$$Cov^Q_0: z \to w \quad \Leftrightarrow \quad
\frac{Q(w)-Q(z)}{w-z}=0.$$ Counting singular points of $f$ now
tells us that $Q$ has one double and two single critical points,
and that therefore up to pre- and post-compositions by M\"obius
transformations $Q$ is the polynomial $Q(z)=z^3-3z$. It follows
that up to conjugacy we may write $f$ in the form
$$z \to w \quad \Leftrightarrow \quad (Jw)^2+(Jw)z+z^2=3.$$
It is easy to see that if we apply a further conjugacy to
transform $J$ to the involution $J(z)=-z$, the equation defining
the correspondence $f$ becomes a member of the family $(2)$. This
completes the proof of Proposition \ref{compatible}.

\subsection{A group-theoretic description of $\Omega$ for a mating}

We shall be pinching unions of arcs in $\Omega(f)$ which are lifts
of simple closed curves in the grand orbit space $\Omega(f)/f$,
where $f$ is one of the matings provided by Theorem 1. With a view
to describing these arcs we examine the structure of $\Omega(f)$
and its relationship with $\Omega(G)$. Our first step will be to
find a {\it Fuchsian} uniformisation for $\Omega(G)/G$.

\medskip
Let $\Gamma$ denote the abstract group $<\sigma,\rho,\tau:
\sigma^2=\rho^3=\tau^2=(\sigma\rho\tau)^2=1>$.

\medskip
Let ${\mathcal F}$ denote the moduli space of conjugacy classes of
faithful discrete co-compact representations of $\Gamma$ in
$PSL_2(\R)$ (recall that a Fuchsian group is said to be {\it
co-compact} if the quotient of the Poincar\'e disc by its action
is compact). An example of a representation of $\Gamma$ which lies
in ${\mathcal F}$ is illustrated in figure 4. Let ${\mathcal F}_0$
denote the path-component of ${\mathcal F}$ containing the
representation illustrated. Thus a faithful discrete
representation of $\Gamma$ lies in ${\mathcal F}_0$ if and only if
there is a fundamental domain $D_\Gamma$ for $\Gamma$ isotopic to
that illustrated in figure 4, with boundary passing through the
fixed points of the corresponding elements of $\Gamma$, in the
same order but with the intervening boundary segments no longer
necessarily geodesic.

\begin{figure}
\begin{center}
\input{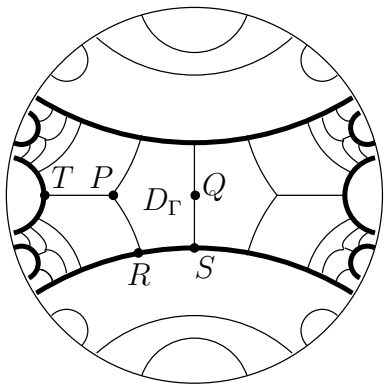}
\caption{A Fuchsian representation of $\Gamma$. Here $P,Q,R$ and
$S$ are the fixed points of $\rho,\sigma,\tau$ and
$\sigma\rho\tau$. The heavy lines indicate the boundary of
$D_{\Gamma_1}$.} \label{Gamma1}
\end{center}
\end{figure}

\medskip
Let $\sigma'=\rho\tau\sigma$. Then $\sigma'$, $\rho$ and $\tau$
together generate $\Gamma=<\sigma,\rho,\tau>$, and satisfy the
same relations. Changing to the new generating set amounts to
applying an (outer) automorphism, which we denote $\beta$, to
$\Gamma$. Let $\psi$ be the automorphism of ${\mathcal F}_0$
induced by composing the representation $\gamma \to PSL_2(\R)$
with $\beta$ and replacing the boundary of $D_\Gamma$ in figure 3
with that given by moving $S$ up to $Q$ (the fixed point of
$\sigma'\rho\tau$), and $Q$ up to $\sigma(S)$ (the fixed point of
$\sigma'$), but keeping $P$ and $R$ unchanged.

\medskip
Recall our description in Section 2.1 of the universal cover
$\widetilde{{\mathcal D}^o}$ of the space ${\mathcal D}^o$ of
conjugacy classes of faithful discrete representations of
$C_2*C_3$ in $PSL_2(\C)$ having connected ordinary set.

\begin{prop}\label{Kleinian-Fuchsian} There is a
homeomorphism $\Psi:\widetilde{{\mathcal D}^o} \to {\mathcal
F}_0$, which carries $t_{1/4}$ to $\psi$ and hence induces
homeomorphisms:

(i) ${\mathcal D}^o \to {\mathcal F}_0/<\psi^2>$;

(ii) ${\mathcal D}^o/\iota \to {\mathcal F}_0/<\psi>$.

\end{prop}

{\bf Proof.} By Lemma \ref{Kleinian-markings}, Section 2.1, a
point of $\widetilde{\mathcal D}^o$ corresponds to an element of
$\tilde{\mathcal S}$, that is to say a sphere equipped with a
complex structure having cone points $P$ of angle $2\pi/3$, and
$Q,R$ and $S$ all of angle $\pi$, together with an isotopy class
of paths $PR$, $RS$ and $SQ$. Obviously it suffices to define a
homeomorphism between $\tilde{\mathcal S}$ and ${\mathcal F}_0$.

\medskip
To do this we uniformise each marked orbifold $\Sigma\in
\tilde{\mathcal S}$ as a quotient of the Poincar\'e disc $\Delta$
by isometries. The marked arcs on $\Sigma$ lift to a union of
arcs, tiling $\Delta$ by translates of a polygon isotopic to that
labelled $D_\Gamma$ in figure 4. The group of covering
transformations of the projection from $\Delta$ to $\Sigma$ is
isomorphic to $\Gamma$ by Poincar\'e's polygon theorem \cite{B}.
Conversely, given a faithful discrete representation of $\Gamma$
lying in ${\mathcal F}_0$, its quotient orbifold $\Sigma$ is an
element of $\tilde{\mathcal S}$. Thus we have a bijection
$\tilde{\mathcal S} \to {\mathcal F}_0$ which, by construction, is
continuous and has a continuous inverse. Since $t_{1/4}$ and
$\psi$ have identical effects on $\Sigma$, our composite
homeomorphism $\Psi:\widetilde{{\mathcal D}^o} \to {\mathcal F}_0$
carries $t_{1/4}$ to $\psi$, and the assertions (i) and (ii) are
immediate corollaries. $\square$

\medskip {\bf Remark.}
The question of finding explicit formulae for bijections between
moduli spaces of representations of Kleinian groups and Fuchsian
groups, such as the bijection provided by Proposition
\ref{Kleinian-Fuchsian}, is in general highly non-trivial, a
classical example being to relate each Schottky group to a
Fuchsian group representing the same surface.

\bigskip Now let $\Gamma_1 \subset \Gamma$ be the subgroup
generated by $\rho\tau$ (which has infinite order), the involution
$\rho^{-1}\tau\rho$, and all involutions of the form
$W\rho^{-1}\tau\rho W^{-1}$, where $W$ runs through those words in
$\sigma$ and $\rho$ which have rightmost letter $\sigma$. Then
$\Gamma_1$ has as fundamental domain the region $D_{\Gamma_1}$
bounded by heavy lines in figure 4. Note that
$D_{\Gamma_1}/\Gamma_1$ is a topological cylinder, the top edge of
the region $D_{\Gamma_1}$ in figure 4 being identified with the
bottom edge, each of the arcs on the left hand edge being folded
in onto an interval, and each of the arcs on the right hand edge
also being folded in onto an interval.

\medskip Suppose $f$ is a $2:2$ holomorphic correspondence which is
a mating, constructed as in Theorem 1, between a faithful discrete
representation of $C_2*C_3$ in $PSL_2(\C)$ having connected
ordinary set and a quadratic map $z \to z^2+c$ having connected
Julia set. Let $\Gamma\subset PSL_2(\R)$ be the Fuchsian
representation associated to it by Lemma \ref{Kleinian-Fuchsian},
and let $\Gamma_1$ be the subgroup of $\Gamma$ defined above.

\begin{prop}\label{groupify} There is a bi-analytic homeomorphism
$$D_{\Gamma_1}/{\Gamma_1}\cong \Delta/{\Gamma_1}
\to \Omega(f)$$ carrying the action of the pair $\{\sigma \rho,
\sigma\rho^{-1}\}$ on $D_{\Gamma_1}/{\Gamma_1}$ to that of the
correspondence $f$ on $\Omega(f)$. \end{prop}

{\bf Proof.} From the construction of the mating $f$ in our
outline proof of Theorem \ref{mating} (in Section 2.2), it is
apparent that $(\Delta,{\Gamma_1})$ uniformises $\Omega(f)$: the
set $D_{\Gamma_1} \cup \rho D_{\Gamma_1} \cup \rho^{-1}
D_{\Gamma_1}$ in figure 4, when quotiented by the boundary
identifications induced by $\Gamma_1$, becomes the annulus $B$ of
figure 3, and the maps $\sigma\rho$ and $\sigma\rho^{-1}$ become
the two `branches' of the correspondence $f$ on $\Omega(f)$.
$\square$

\begin{cor}\label{iota} A mating between $q_c$ and $r\in {\mathcal D}^o$
constructed by the method of Theorem \ref{mating} is canonically
isomorphic to a mating between $q_c$ and $\iota(r)$.
\end{cor}

{\bf Proof.} The outer automorphism $\beta$ defined by replacing
the generator $\sigma$ of $\Gamma$ by $\sigma'=\rho\tau\sigma$
stabilises $\Gamma_1$, and the correspondence induced by
$\{\sigma\rho,\sigma\rho^{-1}\}$ on $\Delta/\Gamma_1$ is the same
as that induced by $\{\sigma'\rho,\sigma'\rho^{-1}\}$, since
$\sigma'\rho=\rho\tau\sigma\rho$ and
$\sigma'\rho^{-1}=\rho\tau\sigma\rho^{-1}$. $\square$

\medskip
{\bf Remarks}

\medskip 1. The idea of regarding $\Omega(f)$ as a quotient of
$\Delta$ by an infinitely generated Fuchsian group is originally
due to Chris Penrose.

\medskip
2. We can recover the action of the Kleinian group
$G=<\sigma,\rho,\chi>$ on $\Omega(G)$ from the action of the
corresponding Fuchsian group $\Gamma=<\sigma,\rho,\tau>$ on
$\Delta$, as follows. Take the polygon
$D_{\Gamma_2}=D_{\Gamma_1}\cup \rho\tau(D_{\Gamma_1})$ formed by
two copies of $D_{\Gamma_1}$, one above the other, identify the
top and bottom edges of this polygon to form a cylinder, then fold
and glue the left-hand edge together and fold and glue the right
hand edge together, to form a sphere. The quotient
$D_{\Gamma_2}/\sim$, which can also be described as an orbit space
$\Delta/\Gamma_2$ for an appropriate infinitely generated subgroup
$\Gamma_2 \subset \Gamma$, is conformally equivalent to
$\Omega(G)$. Indeed $\Gamma_2\cong \pi_1(\Omega(G))$, and the
projection $\Delta \to \Delta/\Gamma_2$ is the universal cover for
$\Omega(G)$. Under the bijection from $D_{\Gamma_2}/\sim$ to
$\Omega(G)$ the ends of $D_{\Gamma_2}$ (the cusps) become the
points of the limit set $\Lambda(G)$ of the action of the Kleinian
group $G$ on $\hat{\C}$.

\section{The pinching deformation}

\subsection{The arcs to be pinched}

To describe the arcs that we shall pinch later, we first fix a
standard faithful discrete representation $r_*$ of $C_2*C_3$
having connected ordinary set, and a path $l$ from a fixed point
$R$ of $\chi\rho$ to a fixed point $S$ of $\chi\sigma$ (so $R$ and
$S$ are as illustrated in figure 2). For convenience we may choose
$r_*$ and $l$ so that the corresponding group $\Gamma$ has the
reflection symmetry in the horizontal axis apparent in figure 4.
Now consider the double cover $\tilde{\Sigma}$ of the orbifold
$\Sigma$ ramified at all four cone points. This is a torus, with a
single cone point $P$ of angle $4\pi/3$, represented by the
central hexagon $D_\Gamma \cup \sigma D_\Gamma$ illustrated in
figure 4, with the top edge identified with the bottom edge, and
the left-hand edge identified with the right-hand edge. While
$\tilde{\Sigma}$ is not itself a quotient of the unit disc
$\Delta$ by a subgroup of $PSL_2(\C)$ (since the cone point is not
of angle $2\pi/n$), nevertheless we may equip ${\tilde{\Sigma}}$
with the metric induced by the restriction of the hyperbolic
metric on $\De$ to the hexagon $D_\Gamma \cup \sigma D_\Gamma$.
The involution  $\si$ (on $\De$) induces an involution
$\tilde\sigma$ on $\tilde{\Sigma}$ such that ${\tilde{\Sigma}}/
\tilde{\sigma}={\Sigma}$.

\begin{lemma}\label{arcs-exist} For each rational number $p/q$
there is a geodesic arc $\delta_{p/q}$ in ${\Sigma}$ which has end
points two of the three cone points of angle $\pi$, which misses
the other cone point of angle $\pi$ and the cone point of angle
$2\pi/3$, and which has lift $\tilde{\delta}_{p/q}$ to
$\tilde{\Sigma}$ a simple closed geodesic of winding number $p/q$.
\end{lemma}

{\bf Proof.} For each such $p/q$ (in lowest terms), there is a
simple closed curve of winding number $p/q$ on the torus
${\tilde{\Sigma}}$, passing through (i) the cone points $Q$ and
$S$ if $q$ is even, (ii) the cone points $Q$ and $R$ if $p$ and
$q$ are both odd, and (iii) the cone points $R$ and $S$ if $p$ is
even. Examples are illustrated in figure 5 for typical cases of
each type.

\begin{figure}
\begin{center}
\input{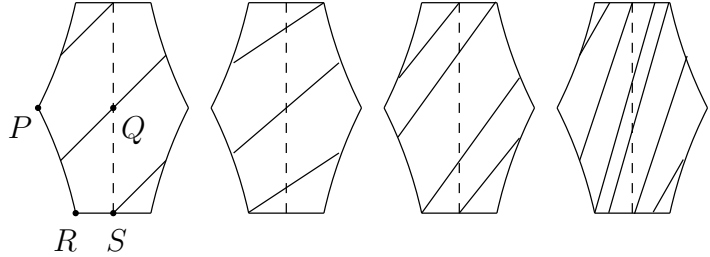}
\caption{The arcs $\tilde{\delta}_{p/q}$ for $p/q=1/2$, $1/3$,
$2/3$ and $4/3$ respectively.} \label{winding-nunmbers}
\end{center}
\end{figure}

\medskip

Note that when we add an even integer to $p/q$ the new
$\delta_{p/q}$ is an arc between the same two cone points on
${\Sigma}$. But when we add an odd integer the roles of $Q$ and
$S$ are interchanged.

\medskip In every case the simple closed curve on
$\tilde{\Sigma}$ can be chosen to be invariant under
$\tilde\sigma$. Since it passes through the lifts of two cone
points, it descends to an arc on ${\Sigma}$ joining these two
points. We define $\delta_{p/q}$ to be a representative of
shortest length in the isotopy class of this arc, relative to its
end points and the other two cone points on ${\Sigma}$.  Note that
there must exist such a minimal length example, as arcs which pass
through one or both of the other cone points have lengths which
are local maxima (since all the cone points have cone angle less
that $2 \pi$). $\square$

\medskip
Let $A_{p/q}$ denote the lift of $\delta_{p/q}$ to the cylinder
$(D_\Gamma \cup \sigma D_\Gamma)/\Gamma_1$ constructed by
identifying the top and bottom of the hexagon. Thus $A_{p/q}$
consists of $q$ arcs each running from one boundary circle of this
cylinder to the other. Consider the union $\Gamma A_{p/q}$ of all
lifts of $\delta_{p/q}$. Recall that $D_{\Gamma_1}/{\Gamma_1}$ is
a cylinder, with ends corresponding to $\partial \Lambda_-$ and
$\partial \Lambda_+$ (by Proposition \ref{groupify}), that the
correspondence $f$ acts on $\partial \Lambda_-$ as a quotient of
the doubling map, and that $f^{-1}$ acts on $\partial \Lambda_+$
as a quotient of the doubling map. For simplicity of description
assume that $\partial \Lambda_-$ is a topological circle and the
action of $f$ on it is that of the doubling map (this is the case
when the quadratic map in the mating corresponds to a value of $c$
in the interior of the main cardioid of the Mandelbrot set):
obvious adaptations are possible for the cases where $\partial
\Lambda_-$ is a proper quotient of the circle.

\begin{figure}
\begin{center}
\input{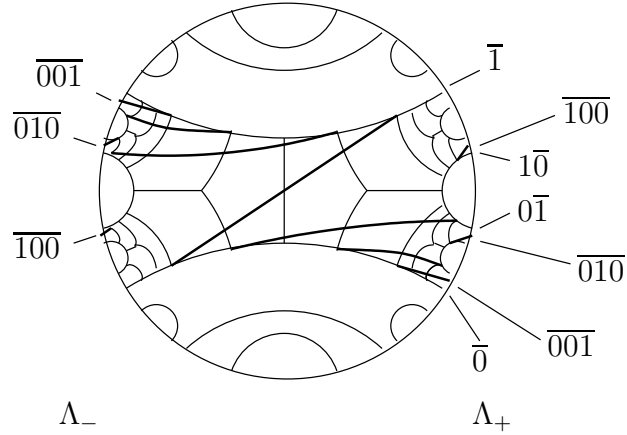}
\caption{The three arcs linking $\Lambda_-$ to $\Lambda_+$ in the case
$p/q=1/3$ (the other images of these arcs under $\Gamma$ are not
shown).}
\label{arcs-on-omega}
\end{center}
\end{figure}

\medskip
If we label the ends of $\partial D_{\Gamma_1}$ by binary
sequences as indicated in figure 6 then the folding
identifications induced by $\Gamma_1$ impose the usual quotient
from the space of binary sequences to the unit circle, carrying
the shift to the doubling map. Thus, under our assumption that
$\partial \Lambda_-$ is the circle, points of $\partial \Lambda_-$
are labelled (figure 6) in such a way that $f^{-1}:\partial
\Lambda_- \to \partial \Lambda_-$ (a $1:2$ correspondence) is
defined by ``right shift and insert $0$ or $1$'' according as the
branch of $f^{-1}$ is $\rho\sigma$ or $\rho^{-1}\sigma$
respectively, and points of $\partial \Lambda_+$ are labelled in
such a way that $f:\partial \Lambda_+ \to \partial \Lambda_+$
(also a $1:2$ correspondence) is defined by ``right shift and
insert $0$ or $1$'' according as the branch of $f$ is $\sigma\rho$
or $\sigma\rho^{-1}$ respectively. We adopt the usual notational
convention that a bar over a symbol (or group of symbols)
indicates the infinite repetition of that symbol (or group of
symbols).

\medskip {\bf Definition} {\it An infinite sequence of $0$'s and
$1$'s is known as {\it Sturmian} if the binary number it
represents on the circle has orbit under the doubling map a
sequence of points arranged in the same order around the circle as
for a rigid rotation.}

\medskip One may assign a rotation number to each Sturmian sequence
$s$, namely the limit as $n$ tends to infinity of the proportion
of the first $n$ digits of $s$ which are $1$'s, or equivalently
the rotation number of the rigid rotation having orbit points in
the same order as those of $s$. Note that such a rotation number
is only defined mod $1$. For each rational $p/q$ (mod $1$) there
is a unique periodic Sturmian orbit of rotation number $p/q$ (this
was observed by Morse and Hedlund, who introduced the notion of
Sturmian sequences). We remark that the points of each periodic
Sturmian orbit ${\mathcal O}$ must be contained in an interval of
length less than $1/2$ on the circle $\R /\Z$, as the doubling map
must preserve the cyclic order of ${\mathcal O}$ (see \cite{BS}
for more about this and other properties of Sturmian sequences).

\medskip
{\bf Examples}

\medskip
The infinite sequences $\overline{01}$, $\overline{001}$ and
$\overline{00101}$ are Sturmian, of rotation numbers $1/2,1/3$ and
$2/5$ respectively.

\begin{prop}\label{landing-points} $(\Gamma A_{p/q}\cap
D_{\Gamma_1})/\Gamma_1$ contains exactly $q$ arcs which join
$\Lambda_-$ to $\Lambda_+$. These land on $\partial\Lambda_-$ at
points of the unique Sturmian orbit of rotation number $p/q$ (mod
$1$) of the $2:1$ map $f:\partial\Lambda_- \to
\partial\Lambda_-$ and at the other end they land on $\partial
\Lambda_+$ at points of the unique Sturmian orbit of $f^{-1}$ of
rotation number $p/q$ (mod $1$).
\end{prop}

{\bf Proof.} The fact that there are exactly $q$ arcs joining
$\Lambda_-$ to $\Lambda_+$ follows at once from the fact that
exactly $q$ arcs in $(\Gamma A_{p/q}\cap D_{\Gamma_1})/\Gamma_1$
cross the equator circle of the central cylinder $(D_\Gamma \cup
\sigma D_\Gamma)/\Gamma_1$ (the vertical line in the central
hexagon in $D_{\Gamma_1}$). The action of the correspondence
$f^{-1}=\{\rho^{-1}\sigma, \rho\sigma\}$ on these arcs is to map
the $j$th arc to the $(j+p)$th arc for each $j$, where the arcs
are counted modulo $q$, from the bottom of the central hexagon
upwards. Thus the action of $f^{-1}$ on the landing point of the
$j$th arc on $\Lambda_+$ is to send it to the landing point of the
$(j+p)$th arc, for each $j$. Similarly $f$ sends the $j$th landing
point on $\Lambda_-$ to the $(j+p)$th. $\square$

\medskip
{\bf Definition of the arc $\gamma_{p/q}$.} {\it For each $p/q$ we
pick as $\gamma_{p/q}$ one of the $q$ components of $(\Gamma
A_{p/q}\cap D_{\Gamma_1})/\Gamma_1$ which cross the equator circle
of the central cylinder and therefore join $\Lambda_-$ to
$\Lambda_+$. For definiteness, when $q$ is odd we take
$\gamma_{p/q}$ to be the component which passes through $R$ (the
fixed point of $\tau$) and when $q$ is even we take it to be the
component which passes through $S$ (the fixed point of
$\sigma\rho\tau$). We remark that in the case $p/q=0$ there is
just one component crossing the vertical symmetry line of the
central hexagon, and it passes through both of these points.}

\medskip
In figure 6 we illustrate $\gamma_{1/3}$, which joins
$\overline{010}\in \Lambda_-$ to $\overline{100}\in \Lambda_+$,
and its two images which also join $\Lambda_-$ to $\Lambda_+$.
These join $\overline{100}\in \Lambda_-$ to $\overline{010}\in
\Lambda_+$, and $\overline{001}\in \Lambda_-$ to
$\overline{001}\in \Lambda_+$ respectively. Arcs
$\gamma_{(3n+1)/3}$ for values of $n$ other than $0$, and their
images, join the same pairs of points in $\Lambda_-$ and
$\Lambda_+$, but wind a different number of times around the
cylinder $D_{\Gamma_1}/\Gamma_1$.

\medskip
For general rational $p/q$ we have the following:

\medskip
{\bf Algorithm} {\it Each point in $\Lambda_-$ represented by a
Sturmian $p/q$ word $\overline{u_1\ldots u_q}$ is joined (by
$\gamma_{p/q}$ or one of its images) to the point in $\Lambda_+$
represented by the Sturmian $p/q$ word
$\overline{u_{q-1}u_{q-2}\ldots u_1u_q}$.}

\medskip
{\bf Proof.} Both $\sigma \rho$ and $\sigma\rho^{-1}$ map the
fixed point $P$ of $\rho$ to $\sigma P$. It follows that $f$ maps
the pair of geodesics landing on $\Lambda_-$ either side of $\bar
1$ to the pair of geodesics landing on $\Lambda_+$ either side of
$\bar 1$ (figure 6). The pair of landing points either side of
$\bar 1$ are represented by the maximum and minimum Sturmian $p/q$
words, $M_{p/q}$ and $m_{p/q}$ respectively, so the arcs landing
at these points of $\Lambda_-$ have their opposite ends at the
points of $\Lambda_+$ represented by $s(m_{p/q})$ and $s(M_{p/q})$
respectively, where $s$ denotes left shift (i.e. `forget the first
digit'). Since it is easily proved from the {\it staircase
algorithm} for Sturmian words \cite{BS} that the minimum word
$m_{p/q}=\overline{v_q\ldots v_1}$ is the reverse of the maximum
word $M_{p/q}=\overline{v_1 \ldots v_q}$, the result follows.
Indeed we may regard the $q$ arcs joining $\Lambda_-$ to
$\Lambda_+$ as indexed by a marked digit in a bi-infinite Sturmian
word, and the action of $f$ and $f^{-1}$ on these arcs as moving
the marker left and right. $\square$

\bigskip {\bf Remarks.}

\medskip
1. Which two of the three cone points on ${\Sigma}$ of cone angle
$\pi$ are the end points of the arc $\delta_{p/q}$ is determined
by the reflection symmetries of the bi-infinite periodic Sturmian
word of rotation number $p/q$ mod $1$. Each such word has
reflection symmetries of exactly two of four possible types:
reflection at a $0$, or at a $1$, or between two adjacent $0$'s or
$1$'s. Which two types occur depends on whether (after reduction
of $p/q$ mod $1$) $p$ is even, $q$ is even, or $p$ and $q$ are
both odd. For example the bi-infinite word generated by
$\overline{00101}$, a case where $p$ is even, has reflection
points between the first two $0$'s and at the third $0$. The
stabiliser of any lift of $\delta_{p/q}$ to $\Delta$ is an
infinite dihedral group, generated by a pair of involutions fixing
adjacent lifts of cone points on the arc, and indeed isomorphic to
the group of symmetries of the bi-infinite periodic Sturmian word.
\medskip

2. The same construction of geodesic arcs crossing the central
hexagon can be followed through for {\it irrational} slope $\nu$
in place of $p/q$. One then obtains a lamination on
$D_{\Gamma_1}/{\Gamma_1}$, with singular leaves passing through
the fixed point of $\rho$ and its translates. In this case the
leaves crossing the hexagon join a Cantor set in $\partial
\Lambda_-$, the unique closed invariant Sturmian set of rotation
number $\nu$ mod $1$, to the analogous Cantor set in $\partial
\Lambda_+$. The algorithm above also applies in this case to tell
us which points are joined to which; we omit details here.

\bigskip It remains to describe the grand orbit of $\gamma_{p/q}$
under the correspondence $f$.

\medskip We start with the special case $p/q=0$. The arc $\gamma_0$
is the lower boundary component of the region $D_{\Gamma_1}$ in
figure 4. Under $f$ this component maps to itself and to the
boundary component of $D_{\Gamma_1}$ which passes through the
point $\sigma(T)$. The grand orbit of $\gamma_0$ under $f$ is the
union of all the boundary components of $D_{\Gamma_1}$, and
quotienting by $f$, or equivalently by $\Gamma_1$, folds all these
components (except the original one) into ``spikes''.

\medskip
We now turn to general $p/q$. From the explicit construction of
matings in Section 2.2 it follows that the branch of $f$ mapping
$\Lambda_-$ to $\Lambda_+$ is defined as follows: given a word $W$
in $0$'s and $1$'s representing a point in $\partial\Lambda_-$ the
$f$-image in $\partial\Lambda_+$ of that point is represented by
the word $\phi(W)$ obtained by changing the parity of the first
digit of $W$. It is now a straightforward computation that when
$q$ is even the set of $q$ arcs joining $\Lambda_-$ to $\Lambda_+$
is mapped two to one by this branch to a set of $q/2$
``concentric'' arcs connecting pairwise the $q$ points of
$\Lambda_+$ obtained by applying the operation $\phi$ to the
Sturmian $p/q$ orbit (i.e. the points of the circle {\it opposite}
to points of the Sturmian orbit). When $q$ is odd, the set of $q$
arcs joining $\Lambda_-$ to $\Lambda_+$ is mapped by this branch
of $f$ to a set of $(q-1)/2$ concentric arcs together with an
innermost spike (figure 7) which lands on $\Lambda_+$ at a single
point, the point opposite to the middle point of the Sturmian
$p/q$ orbit. This spike arises from the fact that for $q$ odd the
geodesic $\gamma_{p/q}$ passes through the fixed point of the
involution $\tau$. Hence its image under the branch of $f$ we are
considering passes through the fixed point of an involution in the
group $\Gamma_1$. This fixed point is on the boundary of
$D_{\Gamma_1}$ (indeed in figure 4 it is the point $\sigma(T)$),
and becomes the end point of a spike in the quotient
$D_{\Gamma_1}/\Gamma_1\cong \Omega(f)$.

\begin{figure}
\begin{center}
\input{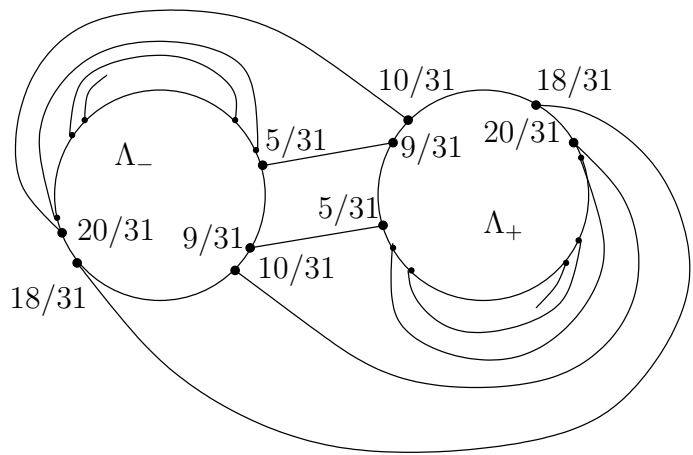}
\caption{The Sturmian orbits of rotation number $2/5$ on
$\Lambda_-$ and $\Lambda_+$, the five arcs joining them, and the
first images of these under the correspondence and its inverse
(subsequent images are not shown).}
\label{images_of_arcs}
\end{center}
\end{figure}

\medskip
Applying $f$ again arbitrarily may times to our ``concentric'' set
of $q/2$ arcs (or $(q-1)/2$ arcs plus a spike, if $q$ is odd), we
obtain smaller and smaller copies around $\partial\Lambda_+$, and
applying $\sigma$ to these copies we obtain similar copies around
$\partial\Lambda_-$, together making up the grand orbit under $f$
of our original set of $q$ arcs.

\subsection{Matings between $q_0$ and circle-packing representations of $C_2*C_3$}

We can now define precisely what we mean by the {\it mating}
between $q_0$ and $r_{p/2q}$ referred to in the statement of \Th
\ref{rat}. After the arcs which make up the grand orbit of
$\gamma_{p/q}$ have been pinched, the intersection $\Lambda_+\cap
\Lambda_-$ is no longer empty, but consists of the $p/q$ Sturmian
orbit of the correspondence on $\partial\Lambda_+$, identified
with the same orbit (in the opposite direction) on
$\partial\Lambda_-$. The set $\Omega$ for the pinched
correspondence has $q$ components whose boundaries meet this
orbit. These form what we call the {\it principal cycle} of
components of $\Omega$. Together with $\Lambda_-\cap\Lambda_+$
itself, they separate the Riemann sphere into two parts, one
containing $\Lambda_-\setminus(\Lambda_+\cap\Lambda_-)$ and the
other containing $\Lambda_+\setminus (\Lambda_+\cap\Lambda_-)$.
The stabilizer (under the iterated pinched correspondence) of each
of the components of the principal cycle is a group, since these
components do not contain ``fold'' points. Moreover it is not hard
to see that this group is isomorphic to $C_2*C_3$.

\medskip
{\bf Definition.} {\it A holomorphic correspondence is said to be
a mating between $r_{p/2q}$ and $q_0$ if it is topologically
conjugate to a correspondence obtained by pinching to a point each
component of the grand orbit of $\gamma_{p/q}$ for a mating
between $r_*$ and $q_0$, and if moreover the action of the
stabiliser of each component of the principal cycle of the
correspondence is conformally conjugate to the action of
$PSL_2(\Z)$ on the upper half-plane.}

\medskip  In a mating between $q_0$ and $r_{p/2q}$, the sets
$\Lambda_+$ and $\Lambda_-$ are no longer copies of $K(q_0)$ (the
unit disc) but are now each homeomorphic to a quotient
$K(q_0)_{p/q}$ of $K(q_0)$ by an equivalence relation $\sim_{p/q}$
on $\partial K(q_0)$ (the unit circle) which may be described as
follows. Let $\omega'_{p/q}$ denote the points of the circle
opposite to points of the Sturmian $p/q$ orbit $\omega_{p/q}$, so
$\omega_{p/q}$ and $\omega'_{p/q}$ are contained in disjoint
intervals. To define the relation $\sim_{p/q}$ we identify the
`outermost' pair of points of $\omega'_{p/q}$, and similarly we
identify the next pair of points from the outside, and so on,
folding the points of $\omega'_{p/q}$ together in pairs. Similarly
we identify in pairs the corresponding inverse images of points of
$\omega'_{p/q}$ under the doubling map, and repeat so that the
relation $\sim_{p/q}$ becomes invariant under this inverse.

\bigskip {\bf Remarks.}

\medskip 1. The justification for describing the construction in
the definition as ``a mating between $q_0$ and $r_{p/2q}$'' is
two-fold. Firstly, both the construction and $r_{p/2q}$ are
obtained by pinching the same simple closed curve $\delta_{p/q}$
on the same orbifold $\Sigma$, and secondly the definition agrees
with our earlier definition for a mating between $q_0$ and the
modular group. However when $p/q\notin {\mathbb Z}$ the most
direct relationship we know of between $\Omega(r_{p/2q})$ and
$\Omega(f)$ for the correspondence pinched along $\gamma_{p/q}$ is
that given by pinching $\delta_{p/q}$ in the Fuchsian picture of
$\Omega(r_*)$, described in Remark 2 following Corollary
\ref{iota} (in Section 2.4).

\medskip
2. Corollary \ref{iota} implies that a mating between $q_0$ and
$r_{p/2q}$ is isomorphic to a mating between $q_0$ and
$r_{(p+q)/2q}$. For example a mating between $q_0$ and $r_{1/2}$
is isomorphic to one between $q_0$ and the modular group. This
example is easily understood directly, since $r_{1/2}$ is the
faithful discrete representation of $C_2*C_3$ for which the limit
set is a single round circle, like $PSL_2(\Z)$, but for which the
generator $\sigma$ of $C_2$ acts by interchanging the two
components of the complement. We remark that $r_{p/2q}$ and
$r_{(p+q)/2q}$ always have the same limit set, since the second
representation is obtained from the first by composing with an
(outer) automorphism of $C_2*C_3$.

\subsection{Invariant collar neighbourhoods of arcs}

For the proofs of \Th \ref{simple} and \Th \ref{rat} we shall need
well-behaved neighbourhoods of our arcs on which to support the
pinching deformations. We define an {\it invariant collar
neighbourhood} of an arc $A$ joining $\Lambda_-$ to $\Lambda_+$ to
be a closed set ${\mathcal N}(A)$ containing $A$, bounded by a
pair of arcs joining the end points of $A$, such that under the
action of $f$ the set ${\mathcal N}(A)$ has stabiliser isomorphic
to the infinite dihedral group, and ${\mathcal N}(A)$ is {\it
precisely invariant} under the action of this stabiliser.
(Strictly speaking, ${\mathcal N}(A)$ is not a topological
neighbourhood of $A$, since the end points of $A$ are on the
boundary of ${\mathcal N}(A)$.)

\begin{lemma}\label{collars_exist}
The arc $\gamma_{p/q}$ has an invariant collar neighbourhood.
\end{lemma}

{\bf Proof.} A collar neighbourhood of each of the $q$ arcs which
join $\Lambda_-$ to $\Lambda_+$ is obtained by lifting any collar
neighbourhood of the $p/q$ geodesic $\delta_{p/q}$ on the orbifold
$\Sigma$. It is immediate from the action of $\sigma\rho$ and
$\sigma\rho^{-1}$ on the lift of such a neighbourhood that its
stabiliser under the action of $f$ is an infinite dihedral group,
generated by the appropriate branch of $f^q$ and by $\sigma$
(which is a branch of $f^{-1}ff^{-1}$) composed with a branch of
whichever $f^r$ maps the $\sigma$ image of the arc back to the
arc. This lifted collar neighbourhood is precisely invariant under
the action of the stabiliser. $\square$

\medskip The small copies of the $q$ arcs have collar
neighbourhoods that are the images of the original collar
neighbourhoods under appropriate branches of forward or backward
iterates of $f$. These images are each either a bijective copy, or
(in the case of a ``spike'') a quotient by an involution, of one
of the original collar neighbourhoods. In the case of the arc
$\gamma_0$, joining the fixed points of the doubling map on
$\partial \Lambda_-$ and $\partial \Lambda_+$, {\it all} the
images are such quotients.

\subsection{A pinching deformation}

Let us consider a correspondence $p$ which represents the mating
of a quadratic polynomial $q$ with a faithful and discrete
representation of $C_2 * C_3$ with connected ordinary set, and let
$f:\La_-\to\La_-$ be the $2:1$-branch of $p$. We fix the curve of
rotation number $p/q$ and consider its lifts $\RRR$ (for red) to
$\cbar$. Thus $\g=\g_{p/q}$ is one of the connected components of
$\RRR$ which joins $\La_-$ to $\La_+$. Let us denote its collar
neighbourhood defined above by $\NNN(\g)$. Then $\stab_p(\NNN(A))$
is isomorphic to the infinite dihedral group. Let $B_-$ and $B_+$
be both components of $\NNN(A)\setminus \g$.

\medskip

We will  first define  an appropriate \qc deformation on a model
strip and then  implement it on the dynamical plane \cite{HT}.


\medskip

Our model space will be a  closed horizontal strip on the upper
half-plane. Choose a collection of numbers $0<  L_y < L_r$ (the
indices $y,r$ are colours  yellow and red respectively), and then
an increasing $C^1$-function $\tau:[0,1[\to [L_r,+\infty[$. Let
$M\subset \R^2$ be the closed subset bounded by
$$([0,1]\times\{0\})\cup (\{0\}\times [0,L_r])\cup (\{1\}\times [0,+\infty[
) \cup ( \{(t,\tau(t)),t\in [0,1[\})\ .$$ Choose $v_t(y)$ so that
$v_t(y)=y$ for $0\le y\le L_y$ and that $(t,y)\mapsto (t,v_t(y))$
is a $C^1$-diffeomorphism from $[0,1]\times [0,L_r]\smallsetminus
\{(1,L_r) \}\to M$.

\begin{figure}
\begin{center}
\input{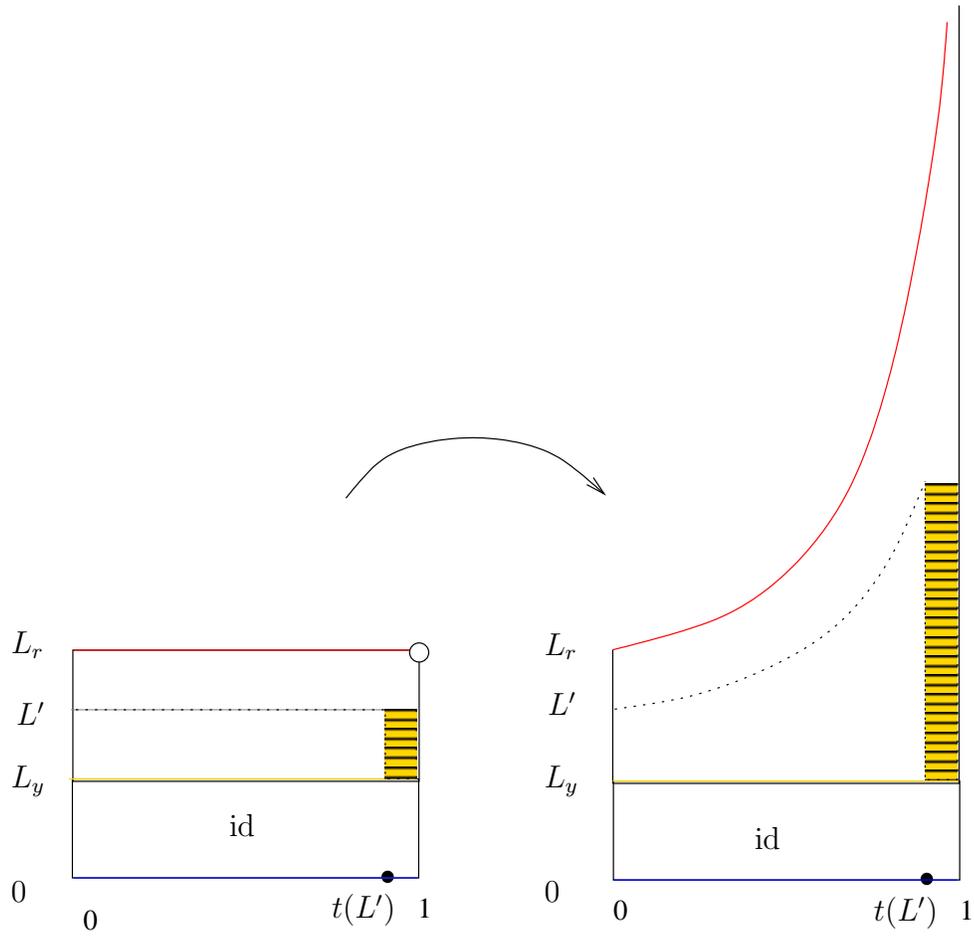}
\caption{The diffeomorphism $(t,y)\mapsto (t, v_t(y))$.}
\label{deform}
\end{center}
\end{figure}

\medskip We also make the following technical assumption: for any
$L'<L_r$, there is  $t(L')\in ]0,1[$ with $t(L')\to 1$ as ${L'\to
L_r}$, such that for any  $(s,y) \in\, ]t(L'),1] \times [0,L']$,
we have $v_s(y)=v_{t(L')}(y)$. Now on the straight strip $\{0\le x
\le L_r\}$,  and for every $t\in [0,1]$ , set
$$\widetilde{P}_t(x+iy)=x+ i\cdot v_t(y)\ .$$ This map satisfies
the following properties\,:

\begin{enumerate}
\item It commutes with the translation by $1$ (and by any other
real number).

\item It is the identity on the sub-strip $\{0\le y\le L_y\}$.

\item The  coefficient of the Beltrami form
$$\left.\frac{\partial \widetilde{P}_t /\partial \bar{z}}{\partial
\widetilde{P}_t /\partial z}\right|_{x+iy}=
\frac{1-\frac{\partial}{\partial
y}v_t(y)}{1+\frac{\partial}{\partial y}v_t(y)}$$ is continuous on
$(t,x+iy)\in [0,1]\times \{0\le y\le L_r\}$, whose norm is locally
uniformly bounded from $1$ if $(t,y)\ne (1,L_r)$ and tends to $1$
as $(t,y)\to (1,L_r)$.

\end{enumerate}

Define conformal maps $\psi_\pm :B_\pm\to \R\times (0,L_r)$ which
map $\g$ to $\R\times\{L_r\}$. For $t\in [0,1[$, set
$\sigma'_t=(\widetilde{P}_t\circ \psi_\pm)^*(\sigma_0)$ to be the
pull-back of the standard complex structure on $B_\pm$. Since the
action is properly discontinuous on $\oo(f)$, we may spread
$\si_t'$ to the whole orbit of $\NNN(\g)$ under the correspondence
$p$. We let $\si_t$ be the extension of this almost complex
structure to the whole Riemann sphere by setting $\si_t=\si_0$ on
the complement. It is a $p$-invariant complex structure. We let
$\YYY$ (for yellow) be the  set of points $z$ such that $\si_t(z)$
is not the standard conformal structure for some $t$.

\medskip

The family of $p$-invariant complex structures $(\si_t)_{t\in
[0,1)}$ defines a pinching deformation supported on $\RRR$. We let
$h_t$ be the \qc map given by the Measurable Riemann Mapping
Theorem applied to $\si_t$ normalised so that $h_t$ fixes both
critical points of $f|_{\La_-}$ and $f^{-1}|_{\La_+}$ and the
point at infinity as well. The correspondence $p_t$ defined by
$p_t(z,w)=p(h_t^{-1}(z),h_t^{-1}(w))$ is holomorphic by
construction, and the family of pairs $(p_t,h_t)_{t\in [0,1)}$
defines a marked pinching deformation.

\section{Convergence of the pinching deformation}

The proofs of both \Th \ref{simple} and \Th \ref{rat} follow
essentially the same lines. We must prove that the pinching
deformation defined in the previous section converges uniformly in
each case, and we must prove that in each case the limit
correspondence has as stabiliser of each of the components of the
principal cycle of $\Omega$ a group conformally equivalent to
$PSL_2(\Z)$. The strategy for proving uniform convergence is
inspired by \cite{H2,HT} where analogous statements are proved for
rational maps and where detailed  proofs can be found.

\medskip

We proceed to prove both theorems simultaneously as far as
possible. We refer to \cite{H2} and \cite{HT} when we can, instead
of repeating the detailed arguments presented in these papers. The
parts of the proofs which differ for the two theorems are
postponed to $4.1$ and $4.2$. In particular we delay the proof of
the key Lemma \ref{nbhd} (stated below). The first step in the
proof of the theorems is to prove that the path of \qc \homeos
$(h_t)$ is equicontinuous. We will apply the following criterion
the proof of which is elementary (cf. Lemma 2.5 in \cite{HT}).

\begin{lemma}{\em\bf (Equicontinuity criterion at a point)}
Let $\HHH=\{h:\DD\to \C\}$ be a family of continuous injective
maps such that $\cup_{h \in\HHH}h(\DD)$ avoids at least 2 points
in $\C$. Let $(U_{n})_{n\ge 0}$ be a nested sequence of disc-like
\nbhds of the origin in the unit disc $\DD$ such that
$A_{n}=\DD\smallsetminus \overline{U_{n}}$ is an annulus. If there
exists a sequence $\eta_n \nearrow +\infty$ such that
$$\forall\,h\in\HHH,\ \forall\,n\ge 0,\ \mbox{\em mod\,}h(A_{n})\,\ge\,
\eta_n,$$
then $\HHH$ is  equicontinuous at the origin.
\end{lemma}

This means that we need to get infinitely many annuli with
controlled moduli. The assumption on the fixed point $\be$ will
give us information on the support of the deformation\,: this will
enable us to prove the following lemma in the respective cases.

\begin{lemma}\label{nbhd}{\em\bf (One good annulus around each Julia point)}
Fix $r>0$.
\begin{itemize}\item[(i)]
For any $x\in \partial\La_-\cup\partial\La_+\smallsetminus\RRR$,
there are  two open neighbourhoods  $N'(x)$ and $N(x)$ of $x$ in $
D(x,\frac{r}{4})$ and $m>0$ such that  ${\rm mod}\,
h_t(N(x)\smallsetminus \overline{N'(x)})\ge m$ for  all $t$.

\item[(ii)] For any $x=\beta_\g \in \RRR\cap (
\partial\La_-\cup\partial\La_+)$, with $\g$ an
$\RRR$-component,  there is a sequence $(t_n)$ in $[0,1)$ tending
to $1$, a nested sequence of annuli $(A_n)_n$ surrounding $\g$,
and a constant $m > 0$ such that  ${\rm mod}\, h_t(A_n)\ge m/n$
for $t\ge t_n$.
\end{itemize}\end{lemma}

Then the weak hyperbolicity condition is used to spread these
annuli at every point and at every scale and therefore to imply
the equicontinuity of $(h_t)$ (cf. the proof of the Proposition
2.3 in \cite{HT} or \S 3 in \cite{H2}). The estimates of the
conformal moduli also enable us to analyse the structure of the
fibres of any limit map and to conclude that its fibres are
exactly the closures of the \ccs of $\RRR$.

\medskip

Any limit $h_1$ satisfies the conclusion of the \th and we may
also extract a convergent sequence $(p_{t_n})$ of the
correspondences to a correspondence $p_1$ (cf. Appendix A in
\cite{HT}).

\medskip

Since the fibre structure is well understood, it follows that if
there are  other limits $(\widehat{h},\widehat{p})$, then
$\widehat{h}\circ h_1^{-1}$ defines a conjugacy which is conformal
off $h_1(\partial\La_-\cup\partial\La_+)$ (cf. Lemma A.2 in
\cite{HT}).

\medskip

Now it can be shown as in \cite{HT} that all the limit
correspondences satisfy the ``weak hyperbolicity'' condition on
the image of $\partial\La_-\cup\partial\La_+$. Since
$\partial\La_-\cup\partial\La_+$ has no interior, a standard
argument of Sullivan implies that the Lebesgue measure of
$h_1(\partial\La_-\cup\partial\La_+)$ is zero (cf.  \Th 4.1
\cite{H1}). Furthermore, the weak hyperbolicity condition on $p_1$
implies that the following rigidity statement holds.

\begin{prop}\label{conjugacy-conformal-off-limits}
Let $p_0$ and $p_1$ be two correspondences which are matings of
weakly hyperbolic \polys with discrete representations of
$C_2*C_3$. If  $p_0$ and $p_1$ are conjugate by a topological
homeomorphism which is conformal off the limit sets, then the
conjugacy is a M\"obius transformation.\end{prop}

The proof of this proposition follows the same lines as
Proposition 6.3 and \Th 0.2 in \cite{H1}. $\square$

\medskip Thus $\widehat{h}\circ h_1^{-1}$ is  a M\"obius
transformation, whence the uniqueness of the limits $(p_t,h_t)$ as
$t$ tends to $1$.

\medskip To complete the proofs of \Th \ref{simple} and \Th
\ref{rat} it now remains only to prove Lemma \ref{nbhd} in both
cases, and to prove that in each case the limit of the family of
pinching deformations corresponds to the mating we are looking
for.

\subsection{The simple case (winding number zero)}

We shall make use of  the statements proved in \cite{HT} for
simple pinchings of rational maps, so we have to show how to get
to that setting.

\medskip

Using McMullen's gluing lemma (Proposition 5.5 in \cite{mc1}), we
may construct a rational map $R$ of degree $2$ which induces a
partition of the sphere $\cbar= K\sqcup \FFF$ where $K$ is the
filled-in Julia set of a quadratic-like map induced by a
restriction of $R$ hybrid-equivalent to $q$, and $\FFF$ is the
basin of attraction of a  fixed point at infinity of multiplier
$1/2$. For the domains of the quadratic-like map, we first choose
a linearising disc $D$ for the point at infinity which contains
the critical value, and set $V=\cbar\setminus \overline{D}$. If
$V'=R^{-1}(V)$, then $R:V'\to V$ is quadratic-like. Furthermore,
we may find a forward-invariant Jordan arc $\kappa$ in $\FFF$
joining the point at infinity with the corresponding $\beta$-fixed
point which only cuts $\partial V$ once, and then transversally.
Let $\widehat{\RRR}$ be the grand orbit  of $\kappa$ for $R$. It
follows that $(\widehat{\RRR}\setminus\kappa)\cap \partial
V=\emty$.

\begin{prop}\label{comp} There is a \qc $\Phi:\cbar\to\cbar$ such that
\begin{itemize}
\item $\Phi(\La_-)=K$ and $\Phi(\RRR)=\widehat{\RRR}$, \item
$\Phi\circ f= R\circ\Phi$ in a \nbhd of $\La_-$, \item
$\overline{\partial}\Phi =0$ a.e. on $\La_-$.
\end{itemize}\end{prop}

{\bf Proof.} We already know that there is a \qc map
$\phi:\cbar\to\cbar$ which fulfills the conclusions of the
Proposition except for the condition on the curves. We let
$U'\subset \subset U$ be simply connected domains such that the
extension $f:U'\to U$ of the branch of the correspondence
$f:\La_-\to\La_-$ is a quadratic-like map hybrid-equivalent to
$q$. It follows from the construction of $f$ that we may assume
that $U$ is a fundamental domain for the involution $J$.
Furthermore, we may also assume that $\phi(U)=V$.

\medskip

We let $\phi_0:U\to V$ be a \qc \homeo isotopic to $\phi$ rel.
$\La_-$ through an isotopy which maps $\partial U$ to
$\partial V$ throughout, and such that $$\phi_0(\g_0\cap
(\overline{U}\setminus U'))= \kappa\cap (\overline{V}\setminus
V')\mbox{ and } R\circ \phi_0|_{\partial U'}=\phi_0\circ
f|_{\partial U'}.$$ This is possible since both sets $U\setminus
\La_-$ and $V\setminus K$ are annuli and since the action of the
maps $f$ and $R$ are 2:1 coverings. Define $(\phi_n)$ inductively
so that $\phi_{n+1}\circ f=R\circ\phi_n$ so that
$\phi_n|_{\La_-}=\phi|_{\La_-}$ and
$\phi_n|_{\overline{U}\setminus U'}=\phi_0|_{\overline{U}\setminus
U'}$. This sequence is a normal family \qc mappings which admits
at least one limit $\Phi:U\to V$. This map satisfies the
conclusion of the proposition. $\square$

\medskip

We now provide a proof of Lemma \ref{nbhd} under the
assumptions of \Th\ref{simple}.

\medskip

{\bf Proof of Lemma \ref{nbhd}.} We first assume that $q$ is not
conjugate  to $z\mapsto z^2 +1/4$. Then by Lemma 2.7 in \cite{HT}
we have the result we seek but for $\widehat{\RRR}$ and the
rational map $R$ in place of $\RRR$ and the correspondence. By
Proposition \ref{comp} this is all we need, except for the case of
the only $\RRR$-component, $\gamma_0$, which is not contained in
the neighbourhood $U$ of $\Lambda_-$. But $\gamma_0$ is a double
cover of any other component $\gamma$ of $\RRR$ by a branch of the
correspondence, and $\gamma_0$ has a neighbourhood which is a
double cover of a disc neighbourhood of $\gamma$, by the same
branch.

\medskip

We now deal with $q(z)=z^2+1/4$. Let us denote by $p$ the mating
of $ q$ with $C_2*C_3$ and let us define $q_0(z)=z^2$, $p_0$ and
$R_0$ the corresponding mating and rational map. We let
$(p_t,\widehat{h}_t)$ be the simple pinching of $p_0$ considered
above, and  $\Phi_0:\La_-(p_0)\to \overline{\DD}$ be given by
Proposition \ref{comp}. It follows from Corollary 3.10 in
\cite{HT} that there is a $\mu$-\homeo, in the sense of David,
$\phi:\C\to\C$, conjugating $p_0|_{\oo(p_0)}$ conformally to
$p|_{\oo(p)}$. Furthermore, a constant $K_0\ge 1$ exists such that
the set of points $z\in \C$ for which the dilatation ratio
$K_\phi(z)$ is at least $K_0$ is contained in the disjoint union
of the orbit of an invariant sector $S\subset int(\La(p_0))$ with
vertex $\be$ (see Lemma 2.1 \cite{H} for details).

\medskip

We claim that the image under $\phi$ of the controlled annuli for
$p_0$ have also controlled moduli. For points outside the red set,
this is because the set where $K_\phi$ is large is contained in
the union of sectors so that the Key lemma in \cite{HT}, which
implies the bounds on the moduli, also holds for these domains.

\medskip

For points in the red set, we must be more precise and use
intermediate results which are established for the proof of Lemma
2.7 in \cite{HT}. We refer to  \S 2.5 in \cite{HT} for the
details. We let $Y$ be the \cc of $\YYY(p_0)$ which contains
$\g_0$. In the proof of the equicontinuity at those points, it is
shown that there  is a sequence $\psi_n: A_n \to
(-C-(n+1),C+(n+1))^2\setminus[-C-n,C+n]^2$ of homeomorphisms,
where $C$ is a fixed positive real number, such that, for $t\ge
t_n$, $\psi_n\circ \widehat{h}_t^{-1}$ is uniformly quasiconformal
off $\YYY\setminus Y$. Moreover, $\psi_n$ maps $\Phi_0(S)\cap A_n$
onto a rectangle $Q_n=[-C-(n+1),-C-n]\times [C_1,C_2]$ for fixed
constants $C_1$ and $C_2$.

\medskip

The bound on the moduli for the cauliflower map $z\mapsto z^2+1/4$
comes from a length-area argument provided by metrics $(\rho_n^t)$
defined as follows. Let $t\ge t_n$; on
$\widehat{h}_t(\YYY(p_0)\setminus Y)$, we let $\rho_n^t=0$ and on
its complement we define
$$\rho_n^t=\frac{1}{|\partial_z \widehat{h}_t\circ\psi_n^{-1}|-|\partial_
{\bar z}\widehat{h}_t\circ\psi_n^{-1}|} \circ (\psi_n\circ
\widehat{h}_t^{-1})\,.$$ This kind of metric is used to  prove the
quasi-invariance of moduli of annuli for \qc maps. This metric
yields the bound $\mod \widehat{h}_t(A_n)\ge m/n$ where $m>0$ is
independent of $n$.

\medskip

Similarly, we let $\widehat{\rho}^t_n=0$ for points in
$h_t\circ\phi(\YYY(p_0)\setminus Y)$ and on the complement, we let
$$\widehat{\rho}^t_n=\frac{\rho_n^t}{|\partial_z \phi_t|-|\partial_{\bar
z} \phi_t|}\circ\phi_t^{-1}\,,$$ where $\phi_t= h_t\circ \phi\circ
\chi_t^{-1}$. It follows from the construction of $\phi$ that
$K_{\phi}\asymp n$ on $Q_n $ (see Lemma 2.1 in \cite{H}), so that
the area of $h_t(\phi(Q_n))$ is at most a multiple of $n$, as the
area of $h_t(\phi(A_n\setminus Q_n))$, for the metric
$\widehat{\rho}^t_n$. Thus, we get $\mod h_t(\phi_0(A_n))\ge c/n$.
Whence we obtain the estimates of the moduli for these points
also. $\square$

\medskip The following proposition now completes the proof of \Th
\ref{simple}.

\begin{prop}\label{simple-mating} Under the assumptions of \Th \ref{simple},
the limit $p_1$ of $(p_t)$ is a mating of $q$ with $PSL_2(\Z)$.
\end{prop}

{\bf Proof.} The limiting correspondence $p_1$ inherits a
compatible involution $J$ from $p_0$, so by Proposition
{\ref{compatible} (Section 1) this correspondence is conjugate to
some member of the family (2), or equivalently to $J \circ
Cov_0^Q$ for $Q(z)=z^3-3z$ and $J$ some (M\"obius) involution. The
proof of the Proposition now follows the same steps as the proof
of Theorem 7.1 in \cite{bf}, which states an analogous result for
the degree $4$ Chebyshev polynomial in place of $Q$. We summarise
the steps but refer the reader to \cite{bf} for technical details.
The topological dynamics of $p_1$ ensure that there exist a
transversal $D_Q$ for $Q$ and a fundamental domain $D_J$ for $J$
such that the complement of the union of the interiors of $D_Q$
and $D_J$ consists precisely of the fixed point $\Lambda_+ \cap
\Lambda_-$. This fixed point is parabolic for $f$ and it follows
from local anaysis that in a neighbourhood the boundaries of $D_Q$
and $D_J$ may be chosen to be smooth curves, tangent to one
another at the fixed point. The set $D_Q\cap D_J$ is a fundamental
domain for the action of $f|_\Omega$, and since $f|_\Omega$ and
$f^{-1}|_\Omega$ have no critical points (only double points) we
know that $f|_\Omega$ is conformally conjugate to
$\{\sigma\rho,\sigma\rho^{-1}\}$ for some Fuchsian representation
of $C_2*C_3$ acting on the open upper half of the complex plane.
To show that this action is indeed that of $PSL(2,\Z)$ it suffices
to show that in the upper half-plane the images of $\partial D_Q$
and $\partial D_J$ converge to the same point on the real axis.
This can be shown to follow from the fact that $\partial D_Q$ and
$\partial D_J$ are smooth curves which meet tangentially (see
\cite{bf}). $\square$

\subsection{Pinching arcs of non-zero rational winding number}

Let $p$ be a correspondence which is a  mating between $z\mapsto
z^2$ and a faithful discrete representation of $C_2*C_3$ in
$PSL_2(\C)$ with connected ordinary set. In this section we prove
Lemma \ref{nbhd} for curves in $\oo(p)$ with non-zero rational
rotation number. The fact that the Julia set is a quasicircle will
be crucial in the proof, which closely follows the argument in \S
3 of \cite{H2}.

\medskip

The first step is to straighten the limit set and the support of
the pinching. Figure \ref{chi} illustrates an example.

\begin{lemma} There is a \qc map $\chi:\cbar\to\cbar$ such that
$\chi(\partial \La_-)=\SS^1$, which satisfies the following
properties\,:
\begin{itemize}
\item $\chi$ is conformal on the interior of $\La_-$\,; \item
$\chi$ conjugates $f$ to $z\mapsto z^2$ in a \nbhd of  the
interior of $\La_-$\,; \item components of $\YYY$ which are
attached at two points $x$ and $y$ to  $\La_-$ are mapped into
rectangles in (log)-polar coordinates with base
$[\chi(x),\chi(y)]$; \item  components $Y$ of $\YYY$ which are
attached at a single point $x$ to $\La_-$ are mapped into sectors
based at $\chi(x)$;
\end{itemize}\end{lemma}

{\bf Proof.} The restriction of $\chi$ to $\La_-$ is given by the
B\"ottcher coordinates of $f$. The extension of $\chi$ to the
outside makes use of a pull-back argument (see pp.\,14-15 in
\cite{H2}). $\square$

\medskip

\begin{figure}
\begin{center}
\input{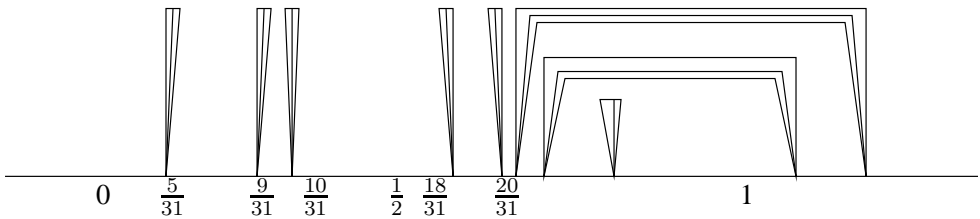}
\caption{Image under $\chi$ of the collars of the first two
generations of the orbit of $\gamma_{p/q}$, in the case $p/q=2/5$
(cf. fig. 7).} \label{chi}
\end{center}
\end{figure}

\medskip

The next step of the proof is to control the moduli of many
annuli. We place ourselves in the coordinates given by $\chi$. As
in \cite{H2}, we may define annuli bounded by rectangles in the
log-polar coordinates which avoid the image of $\YYY$ under
$\chi$.

\medskip
As in the case of simple pinchings, there is no problem with the
curves which link both components of $\La$, because they cover
other components which do not. This enables us to prove Lemma
\ref{nbhd} (cf. Proposition 3.3 and 3.4 in \cite{H2}).

\medskip
Finally, the following proposition completes the proof of \Th
\ref{rat}.

\begin{prop}\label{rational-mating} Under the assumptions of
\Th \ref{rat}, the limit $p_1$ of $(p_t)$ is a mating of $z\mapsto
z^2$ with the circle-packing representation $r_{p/2q}$ of
$C_2*C_3$.
\end{prop}

{\bf Proof.} As in the proof of Proposition \ref{simple-mating}
the limiting correspondence $p_1$ is necessarily conjugate to some
member of the family (2), or equivalently to $J \circ Cov_0^Q$ for
$Q(z)=z^3-3z$ and $J$ some (M\"obius) involution. Once again we
can now follow the same steps as in the proof of Theorem 7.1 in
\cite{bf}. Transversals $D_Q$ and $D_J$ can be chosen this time
such that the complement of the union of their interiors consists
precisely of the period $q$ parabolic orbit $\Lambda_+ \cap
\Lambda_-$, and such that in a neighbourhood of any point of this
orbit the boundaries of these transversals are smooth curves,
tangent to one another at the orbit point. From the fact that
$\Omega$ is now a countable union of topological discs and our
knowledge of the topological dynamics of $f$ (using convergence of
the pinching deformation) we know that $f|_\Omega$ and
$f^{-1}|_\Omega$ have no critical points (only double points) and
that for any component of $\Omega$ which meets the period $q$
orbit $\Lambda_+\cap\Lambda_-$ the iterated branches of $f$ which
stabilise the component are conformally conjugate to the elements
of the group generated by $\{\sigma\rho,\sigma\rho^{-1}\}$ for
some Fuchsian representation of $C_2*C_3$ acting on the open upper
half of the complex plane. As in the proof of Proposition
\ref{simple-mating} the properties of the boundaries of $D_Q$ and
$D_J$ again ensure that this representation is indeed conformally
conjugate to $PSL_2(\Z)$. $\square$

\end{document}